\documentclass[10pt]{amsart}
\usepackage[arrow,matrix,curve]{xy}
\usepackage{euscript,amssymb,amsbsy}
\usepackage{graphicx}

\newtheorem{stuff}{Stuff}[section]
\newtheorem{theorem}[stuff]{\sl Theorem}
\newtheorem{proposition}[stuff]{\sl Proposition}
\newtheorem{lemma}[stuff]{\sl Lemma}
\newtheorem{corollary}[stuff]{\sl Corollary}

\newenvironment{definition}{%
\vskip1ex\refstepcounter{stuff}\trivlist \itemindent 0pt
\item[\hskip\labelsep\sl Definition \thestuff.]%
\ignorespaces}{\endtrivlist\vskip1ex}%
\newenvironment{remark}{%
\vskip1ex\refstepcounter{stuff}\trivlist \itemindent 0pt
\item[\hskip\labelsep\sl Remark \thestuff.]%
\ignorespaces}{\endtrivlist\vskip1ex}%

\let\rar\rightarrow
\let\lar\longrightarrow

\let\dar\downarrow

\let\hra\hookrightarrow
\let\mt\mapsto
\let\lmt\longmapsto
\let\xrar\xrightarrow

\let\euf\EuScript 
\let\cal\mathcal
\let\mbb\mathbb

\DeclareFontFamily{OT1}{rsfs}{}
\DeclareFontShape{OT1}{rsfs}{n}{it}{<->rsfs10}{}
\DeclareMathAlphabet{\crl}{OT1}{rsfs}{n}{it}

\long\def\InsertFig#1 #2 #3 #4 #5\EndFig{\scalebox{#1}{\hskip #2 mm
$\vbox to #3mm{\vfil\includegraphics{#4}}#5$}}
\long\def\LabelTeX#1 #2 #3\ELTX{\rlap{\kern#1mm\raise#2mm\hbox{#3}}}

\let\unl\underline

\let\tld\tilde

\let\nit\noindent

\let\disp\displaystyle
\let\srel\stackrel
\let\sm\setminus

\let\veps\varepsilon

\numberwithin{equation}{section}

\newcommand\Hom{\mathop{\sf Hom}\nolimits}
\def\Aut{\mathop{\sf Aut}\nolimits}
\def\Ker{\mathop{\sf Ker}\nolimits}
\def\Sym{\mathop{\sf Sym}\nolimits}
\let\si\sigma
\let\Si\Sigma
\newcommand\invq{{\slash\kern-0.65ex\slash}}
\newcommand\pr{{\rm pr}}

\newcommand\uz{\unl z}
\newcommand\hg{H^*}
\newcommand\LR{\langle{\rm LR}\rangle}
\newcommand\SR{\langle{\rm SR}\rangle}
\newcommand\cLR{\langle{\cal{LR}}\rangle}
\newcommand\cSR{\langle{\cal{SR}}\rangle}
\newcommand\BG{{\rm BG}}

\newcommand\bfd{{\hbox{$\boldsymbol{\cdot}$}}}
\newcommand\csi{{{({\mbb C}^*)}^{\Si(1)}}}
\newcommand\cxi{{{({\mbb C}^*)}^\Xi}}
\newcommand\tcl{{\kern1ex\tld{\kern-1ex\crl L}}}
\newcommand\tcu{{\kern1ex\tld{\kern-1ex\crl U}}}
\newcommand\tcx{{\kern1ex\tld{\kern-1ex\crl X}}}
\begin{document}

\title{Families of toric varieties}
\author{Mihai Halic}
\address{Institut f\"ur Mathematik, Universit\"at Z\"urich, 
Winterthurerstrasse 190, CH-8057 Z\"urich, Switzerland}
\email{halic@math.unizh.ch}
\maketitle
\markboth{\sc FAMILIES OF TORIC VARIETIES}{MIHAI HALIC}

The purpose of this note is to present a construction for 
families of complete toric varieties, which extends the 
standard one of a single variety. The construction we are 
going to describe is very much in the spirit of \cite{cox}, 
where toric varieties are described as quotients of open 
subsets in affine spaces: here we simply replace the affine 
space by the total space of a vector bundle over an 
arbitrary base variety.  

Such families naturally occur in the study of moduli 
spaces of curves on toric varieties (see \cite{ha}) or, 
alternatively, as the space of solutions of certain 
vortex-type equations on a Riemann surface (see \cite{ot}). 
Integration on these moduli spaces do furnish invariants, 
in this case for the toric variety we start with. It is 
therefore important to be able to perform explicit 
computations in cohomology, and this in turn requires a 
good description of the cohomology ring of the moduli space.

With this motivation in mind, we will construct and will 
study families of toric varieties for their own sake. 
The structure of the paper is a very linear one: as usual, 
the first section is devoted to set up the stage. 
We review the construction of toric varieties as quotients 
of open subsets in affine spaces, but it is  our intension 
to emphasize the coordinate free features of it. The 
coordinate free treatement will eventually enable us 
to extend the construction to the relative setting. 

This is done in the second section: we prove that is possible 
to `paste' {\em complete} toric varieties together, and to 
obtain a relative version of the standard construction of 
toric varieties. The objects we get in this way are categorical 
quotients of Zariski open subsets in the total space of a 
vector bundle over an arbitrary base, and they are natural 
generalizations of the projectivized bundle of a vector bundle.

The next step is to compute the cohomology ring of families 
of complete and simplicial toric varieties. We make use of 
the computation in \cite{da} of the cohomology ring of complete 
and simplicial toric varieties and we find that the structure 
of the cohomology ring in the relative setting is similar to 
the standard case. 

We conclude the article with an overwiew of the situation 
appearing in \cite{ha,ot} in the context of the gauged sigma 
models for toric varieties, which was our primary motivation 
for this study, and we place it into the more formal frame 
of the present paper.\medskip 

{\small\it Acknowledgements}\quad I wish to thank prof. Ch.\ Okonek 
for encouraging me to develop this subject, and for his critical 
remarks.


\section{Some preliminary results}{\label{sct:rmk}}

As we have already said in the introduction, we start with 
a remark concerning the primitive families of complete fans. 
Let $\Si\subset N_{\mbb R}$, with $N\cong\mbb Z^n$, be a 
complete fan and denote by $\Xi$ the set of generators of 
the $1$-dimensional cones of $\Si$: for each $\xi\in\Xi$, 
${\mbb R}_{\geq 0}\xi\in\Si(1)$ and $\xi$ is the generator 
of the semi-group ${\mbb R}_{\geq 0}\xi\cap N$; moreover, 
we denote by $(f_\xi)_{\xi\in\Xi}$ the canonical basis of 
${\mbb Z}^\Xi$. Putting $M:=N^\vee$, we obtain the standard 
short exact sequence 

\begin{align}{\label{sqn:exact}}
\begin{array}{l}
0\lar M\srel{a}{\lar} {\mbb Z}^\Xi\srel{c}{\lar} A_\Xi\lar 0,\quad 
\text{where }A_\Xi:={\rm Coker}(a),\\ 
\hskip 6.5ex\disp m\lmt \sum_{\xi\in\Xi}\langle m,\xi\rangle f_\xi.
\end{array}
\end{align}

\begin{definition}
We say that $\xi',\xi''\in\Xi$ are equivalent, and we write 
$\xi'\sim\xi''$, if $c(f_{\xi'})=c(f_{\xi''})$.
\end{definition}

Our first goal in this section is to prove the

\begin{proposition}{\label{prop:prim}}
Let $\Si$ be a complete fan and consider $\xi',\xi''\in\Xi$. 
If $\pi\subset\Xi$ is a primitive family, $\xi'\in\pi$ and 
$\xi'\sim\xi''$, then $\xi''\in\pi$ too.
\end{proposition}

Before starting the proof, let us notice that in the proposition 
the restriction on $\Si$ to be complete is essential. Consider for 
instance the regular fan
$$
\begin{array}{ccl}
\InsertFig
1 14 0 {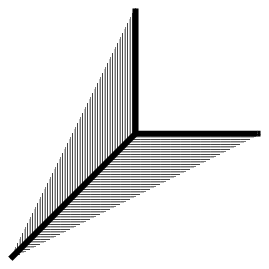} 
{\LabelTeX  27 12 $\xi_1$\ELTX}
{\LabelTeX  14 25 $\xi_2$\ELTX}
{\LabelTeX -15 -4 $-(\xi_1+\xi_2)=\xi_3$\ELTX}
\EndFig
&\hskip3cm &
\raise1cm\hbox{
\begin{minipage}[b]{39ex}
In this case the primitive family is\break 
$\pi=\{\xi_1,\xi_2\}$, 
but $c(f_{\xi_3})=c(f_{\xi_1})=c(f_{\xi_2})$.
\end{minipage}}
\end{array}
$$

\begin{proof}
Let $R={\mbb C}[Z_\rho;\rho\in\Xi]$ be the homogeneous 
coordinate ring of the toric variety $X_\Si$ defined by $\Si$. 
The group $G_\Xi:={\sf Hom}_{\mbb Z}(A_\Xi,{\mbb C}^*)$ 
acts on $R$ and the weights of the action define a graduation 
of $R$ (see \cite{cox}, section 1). Let us denote by
$$
R=\bigoplus_{\beta\in A_\Xi}R_\beta
$$
the graded pieces of $R$, and by $\Aut_g(R)$ the group of graded 
automorphisms of $R$. We know by \cite{cox}, theorem 2.1, that 
$X_\Si$ is canonically isomorphic to the categorical quotient 
$U(\Si)/G_\Xi$, where $U(\Si)$ is the complement in ${\mbb C}^\Xi$ 
of the `bad locus' 
$$
Z_\Si=\bigcup_{\text{$\pi$ primitive family}}\kern-2ex{\mbb A}(\pi),
\quad 
{\rm with}\; {\mbb A}(\pi):=\{ z_\rho=0,\, \forall\rho\in\pi\}.
$$
We have the following informations concerning $\Aut_g(R)$:
\smallskip

\nit\unl{\sl Result}\;\cite[prop. 4.3, 4.6]{cox}\; 
{(i) \it $\Aut_g(R)$ is a connected algebraic group;

\nit\rm (ii) \it Each $\phi\in\Aut_g(R)$ defines an automorphism 
$\tld\phi$ of ${\mbb C}^\Xi$ having the property that 
$\tld\phi(Z_\Si)=Z_\Si$.
}\smallskip

Actually, the second statement is proved in \cite{cox} only for 
complete and simplicial fans. However, in the appendix A we 
show that the proof given in {\it loc. cit.} extends to complete 
fans too.

Let us notice that since 
$\disp Z_\Si=\kern-2ex\bigcup_{\text{$\pi$ primitive}}
\kern-2ex{\mbb A}(\pi)$ 
is the decomposition of $Z_\Si$ into its irreducible components 
and since the group $\Aut_g(R)$ is connected and preserves 
$Z_\Si$, each component ${\mbb A}(\pi)$ must be preserved under 
the $\Aut_g(R)$-action. 

Let us consider now $\xi', \xi''\in \Xi$ with $\xi'\sim\xi''$ 
and assume that $\pi_0$ is a primitive family with 
$\xi'\in\pi_0$; we want to prove that $\xi''$ still belongs 
to $\pi_0$. Since $c(f_{\xi'})=c(f_{\xi''})$, the assignment
$$
\begin{array}{l}
\phi(Z_\xi):=Z_\xi,\quad\forall\xi\in\Xi\sm\{\xi',\xi''\},\\ 
\phi(Z_{\xi'}):=Z_{\xi''}
\quad\text{and}\quad
\phi(Z_{\xi''}):=Z_{\xi'}
\end{array}
$$
extends to a graded automorphism of $R$, that is it defines 
an element of $\Aut_g(R)$. The only effect of the corresponding 
automorphism $\tld\phi$ of ${\mbb C}^\Xi$ is to exchange the 
coordinates $z_{\xi'}$ and $z_{\xi''}$. Let us assume now that 
$\xi''\not\in\pi_0$. Then, on one hand, 
$$
\tld\phi\bigl({\mbb A}(\pi_0)\bigr)
=\tld\phi\bigl(\{z_\xi=0\mid\xi\in\pi_0\}\bigr)
=\left\{\left.
\begin{array}{l}
z_\xi=0\\ z_{\xi''}=0
\end{array}\right|
\xi\in\pi_0\sm\{\xi'\}
\right\},
$$
and, on the other hand, ${\mbb A}(\pi_0)$ is invariant under the 
group action 
$$
\tld\phi\bigl({\mbb A}(\pi_0)\bigr)
={\mbb A}(\pi_0)
=\left\{\left.
\begin{array}{l}
z_\xi=0\\ z_{\xi'}=0
\end{array}\right|
\xi\in\pi_0\sm\{\xi'\}
\right\}.
$$
Clearly, the two equalities are contradictory, so it must 
be that $\xi''\in\pi_0$.
\end{proof} 

\begin{definition}{\label{def:incl}}
For a primitive family $\pi$, we say that an element 
$\beta\in A_\Xi$ is included in $\pi$, and we write 
$\beta\subset\pi$, if there is a $\xi'\in \Xi$ such 
that $\xi'\in\pi$ and $c(f_{\xi'})=\beta$. 
\end{definition}

The previous proposition implies that $\beta\subset\pi$ if 
and only if $\{\xi\mid f_\xi\in c^{-1}(\beta)\}\subset\pi$. 
In other words, the definition above does make sense.
\smallskip

We pass now to the second part of this section, where we recall 
the construction of toric varieties as quotients of open subsets 
in affine spaces. Here we will work in a more general setting 
than before: $N\cong{\mbb Z}^n$ is going to be a lattice again, 
but this time $\Xi\subset N$ is a set which satisfies the 
following properties: (i) $0\not\in\Xi$, 
(ii) $\sum_{\xi\in\Xi}{\mbb R}_{\geq 0}\xi=N_{\mbb R}$, 
(iii) every $\xi\in\Xi$ is the generator over 
${\mbb Z}_{\geq 0}$ of the semi-group 
${\mbb R}_{\geq 0}\xi\cap N$. 

We are interested in complete fans $\Si\subset N_{\mbb R}$ 
having the property that the generators of $\Si(1)$ form a 
subset of $\Xi$. Let us fix such a fan $\Si$, and denote by 
$J:=\Xi\sm \Si(1)$ the complement of $\Si(1)$ in $\Xi$. In 
analogy with \eqref{sqn:exact}, we have this time the 
commutative diagram 
\begin{align}{\label{cd1}}
\xymatrix{
 & & 0\ar[d] & 0\ar[d] & 
\\ 
  & &{\mbb Z}^J\ar[r]^-\cong\ar[d]& A_J\ar[d]& 
\\ 
0\ar[r]& M\ar[r]^-{a_\Xi}\ar@{=}[d]& 
{\mbb Z}^\Xi\ar[r]^-{c_\Xi}\ar[d]&
A_\Xi\ar[r]\ar[d]^-b& 0 
\\ 
0\ar[r]& M\ar[r]^-{a_\Si}& 
{\mbb Z}^{\Si(1)}\ar[r]^-{c_\Si}\ar[d]&
A_\Si\ar[r]\ar[d]& 0 
\\ 
 & & 0 & 0 & 
}
\end{align}
and by dualizing we obtain the commutative diagram of groups 
\begin{align}{\label{cd2}}
\xymatrix{
 & 1\ar[d] & 1\ar[d] & & 
\\
1\ar[r] & G_\Si\ar[r]\ar[d]& 
{({\mbb C}^*)}^{\Si(1)}\ar[r]\ar[d]&
T_N\ar[r]\ar@{=}[d]& 1 
\\ 
1\ar[r] & G_\Xi\ar[r]\ar[d]& 
{({\mbb C}^*)}^\Xi\ar[r]\ar[d]&
T_N\ar[r]& 1
\\ 
 & G_J\ar[r]^-\cong\ar[d] & {({\mbb C}^*)}^J\ar[d] & & 
\\
 & 1 & 1 & & 
} 
\end{align}

We are interested in describing the decompositions of 
${\mbb C}^\Xi$ into the weight spaces corresponding 
to the $G_\Xi$ and $G_\Si$-actions. The $G_\Xi$-action 
determines the decomposition
$$
{\mbb C}^\Xi=\bigoplus_{\xi\in\Xi}{\mbb C}f_\xi
=\bigoplus_{\alpha\in A_\Xi}
\biggl(
\bigoplus_{f_\xi\in c_\Xi^{-1}(\alpha)}{\mbb C}f_\xi
\biggr)
=:\bigoplus_{\alpha\in A_\Xi}V_\alpha,
$$
while for the $G_\Si$-action we have
$$
{\mbb C}^\Xi=\bigoplus_{\beta\in A_\Si}
\biggl(
\bigoplus_{\alpha\in b^{-1}(\beta)}V_\alpha
\biggr)
=:\bigoplus_{\beta\in A_\Si}W_\beta.
$$
We can see from the second diagram that induced $G_\Si$-action 
on the ${\mbb C}^J$ summand of ${\mbb C}^\Xi={\mbb C}^{\Si(1)}
\oplus{\mbb C}^J$ is trivial. The following lemma says that 
${\mbb C}^J$ can be intrinsically characterized as the trivial 
isotypic component of ${\mbb C}^\Xi$ corresponding to the 
$G_\Si$-action.

\begin{lemma}{\label{lm:triv-subsp}}
{\rm (i)} ${\mbb C}^{\Si(1)}=\bigoplus_{\beta\in 
A_\Si\sm\{0\}}W_\beta$ 
and ${\mbb C}^J=W_0=\bigoplus_{j\in\ker(b)}V_j$;

\nit {\rm (ii)}  Every non-zero weight space $V_j$, with 
$j\in\ker(b)$, has dimension one.
\end{lemma}

\begin{proof}
For the first part it is sufficient to show that there are no 
linear subspaces of ${\mbb C}^{\Si(1)}$ on which $G_\Si$ acts 
trivially. If there was such a linear subspace, acted on 
trivially by $G_\Si$, we deduced the existence of a 
$\xi_0\in\Si(1)$ with the property that $c_\Si (f_{\xi_0})=0$. 
This would imply the existence of an element $m\in M$ such that 
$\langle m,\xi_0\rangle=1$ and $\langle m,\xi\rangle=0$ for all 
$\xi\in\Si(1)\sm\{\xi_0\}$. But this contradicts the assumption 
that $\Si$ is complete.

For the second claim, let us assume the contrary, that there 
exists an $j\in\ker(b)$ such that $\dim V_j\geq 2$. 
Then, from the first part of this lemma, follows that there 
are $\eta', \eta''\in J$ such that 
$c_\Xi(f_{\eta'})=c_\Xi(f_{\eta''})=j$. This in turn implies 
the existence of $m\in M\sm\{0\}$ such that 
$$
\langle m,\eta'\rangle=-\langle m,\eta''\rangle=1 \text{ and }
\langle m,\xi\rangle=0,\;\forall\ \xi\in\Xi\sm\{\eta',\eta''\}.
$$
In particular, $\langle m,\xi\rangle=0$ for all $\xi\in\Si(1)$, 
and this contradicts that $\Si$ is a complete fan.
\end{proof}

Using this language, we can describe more invariantly the 
`bad set' in ${\mbb C}^{\Si(1)}$: 
\begin{align}{\label{eq:equiv}}
\begin{array}{ll}
\disp 
Z_\Si=\bigcup_{\pi\text{ primitive}}\kern-2ex{\mbb A}(\pi), 
\quad\text{with}\quad 
{\mbb A}(\pi)
&
=\{\uz_\alpha:=(z_\xi)_{f_\xi\in c_\Si^{-1}(\alpha)}=0, 
\forall\, \alpha\subset\pi\}
\\[1.5ex]
&
=\Ker\bigl[\pr:{\mbb C}^{\Si(1)}\rar 
\oplus_{\beta\subset\pi} W_\beta\bigr]
\end{array}
\end{align}

We have already recalled from \cite{cox}, theorem 2.1, that 
when $\Xi=\Si(1)$, the toric variety $X_\Si$ can be described 
as the categorical quotient $U(\Si)/G_\Si$. It is obtained by 
glueing together the affine pieces $X_\si=U_\si/G_\Si$, 
$\si\in\Si$, where $U_\si:=\{\hat Z_\si\neq 0\}\subset 
{\mbb C}^{\Si(1)}$ and $\hat Z_\si:=\prod_{\xi\in\Si(1)
\sm\si(1)}Z_\xi$. 

We want to generalize this situation to our context where the 
fan $\Si$ has the property that $\Si(1)\subseteq\Xi$ only. 
Let us consider the projection 
$q:{\mbb C}^\Xi\rar {\mbb C}^{\Si(1)}$ corresponding to the 
decomposition ${\mbb C}^\Xi={\mbb C}^{\Si(1)}\oplus{\mbb C}^J$, 
and $W_0^\times:={({\mbb C}^*)}^J\subset{\mbb C}^J$.

\begin{proposition}{\label{prop:inv-quot}}
The toric variety $X_\Si$ is isomorphic to the categorical 
quotient $\bigl(U(\Si)\times W_0^\times\bigr)/G_\Xi$. If $\Si$ 
is simplicial, this quotient is geometric.
\end{proposition}

\begin{proof}
What we shall actually prove is that the categorical quotient 
$\bigl(U(\Si)\times W_0^\times\bigr)/G_\Xi$ is isomorphic to 
$U(\Si)/G_\Si$, 
and the conclusion will follow then from the result of D.~Cox.

The open subset $U(\Si)\times W_0^\times\subset{\mbb C}^\Xi$ is 
the union of the affine pieces $U_\si\times W_0^\times$, and it 
suffices to check that the quotients 
$\bigl(U_\si\times W_0^\times\bigr)/G_\Xi$ 
and $U_\si/G_\Si$ are isomorphic, and that these isomorphisms 
can be glued together on the overlaps. For $\si\in\Si$ a cone, 
the ring of functions on $U_\si\times W_0^\times$ is 
$$
{\mbb C}[U_\si\times W_0^\times]
={\mbb C}[U_\si]\otimes_{\mbb C}{\mbb C}[W_0^\times]
=R_{\hat Z_\si}\otimes_{\mbb C}
{\mbb C}[Z_j,Z_j^{-1}\mid j\in J]_{Z_J},
$$
where the lower indices denote the localizations at 
$\hat Z_\si=\prod_{\xi\in\Si(1)\sm\si(1)}Z_\xi$ and 
$Z_J:=\prod_{j\in J}Z_j$ respectively. Using the 
diagram \eqref{cd1}, we deduce that
$$
{\mbb C}[U_\si\times W_0^\times]^{G_\Xi}
=
\bigoplus_{m\in \si^\vee\cap M}
\mbb C\ {\cdot} \prod_{\rho\in\Xi} Z_\rho^{\langle m,\rho\rangle}.
$$
and that 
$$
{\mbb C}[U_\si]^{G_\Si}
=
\bigoplus_{m\in \si^\vee\cap M}
\mbb C\ {\cdot} \prod_{\xi\in\Si(1)} Z_\xi^{\langle m,\xi\rangle}.
$$
We define the ring homomorphism
$$
{\mbb C}[U_\si\times W_0^\times]^{G_\Xi}\lar 
{\mbb C}[U_\si]^{G_\Si}\text{ by } 
P=P(Z_\xi,Z_j)\lmt Q_P:=P(Z_\xi,1).
$$
It is both injective and surjective because $\Si$ is complete. 
The isomorphisms clearly patch over the overlaps, and therefore 
we deduce the isomorphism of the two categorical quotients.

Let us assume now that $\Si$ is simplicial. Then we know that 
$U(\Si)/G_\Si$ is a geometric quotient. Since for any $\si\in\Si$,
$$
{\mbb C}[U_\si\times W_0^\times]^{G_\Xi}
=
{\biggl[
{\mbb C}[U_\si\times W_0^\times]^{G_\Si}
\biggr]}^{G_\Xi/G_\Si=G_J}
=
{\biggl[
{\mbb C}[U_\si]^{G_\Si}\times_{\mbb C} {\mbb C}[W_0^\times]
\biggr]}^{G_J},
$$
we can decompose the $G_\Xi$-quotient into a sequence of 
two succesive quotients
$$
\xymatrix@C=-5ex{
U(\Si)\times W_0^\times\ar[rr]^-{G_\Si}\ar[dr]_-{G_\Xi}
&& 
X_\Si\times W_0^\times\ar[dl]^-{G_J}
\\ 
&\bigl(U(\Si)\times W_0^\times\bigr)/G_\Xi.&
}
$$
We know that the $G_\Si$-quotient is geometric, and we notice 
that $G_J$ acts freely and transitively on 
$W_0^\times\subset{\mbb C}^J\subset{\mbb C}^\Xi$. Therefore 
the second quotient is still geometric. Globally, we deduce 
that the $G_\Xi$-quotient is geometric.
\end{proof}

We must notice that the isomorphism between the quotients 
$\bigl(U(\Si)\times W_0^\times\bigr)/G_\Xi$ and $U(\Si)/G_\Si$ 
is {\em not} natural, and this fact will have as a consequence 
that the two quotients will give rise, in the relative setting, 
to distinct families of toric varieties.

Another point that we wish to discuss is related to the sheaves 
associated to the invariant divisors on a toric variety. It is 
well-known that given a fan $\Si$, to each element $\xi\in\Si(1)$ 
we can associate an irreducible Weil divisor $D_\xi\subset X_\Si$, 
invariant under the action of the big torus of $X_\Si$, and any 
such divisor defines the rank one, coherent sheaf 
${\cal O}(D_\xi)\rar X_\Si$. The shortcomming of this construction 
is that it depends on coordinates, though the sheaf 
${\cal O}(D_\xi)$ itself does not: for $\xi'\sim \xi''$, the 
corresponding sheaves are isomorphic. Our next goal is to 
describe ${\cal O}(D_\xi)$ in a coordinate free fashion, which 
will eventually allow us to define their analogues in the 
relative setting.

\begin{proposition}{\label{prop:coord-f}}
Let $\Xi\subset N$ be a finite set obeying the conditions 
(i)--(iii), and let $\Si\subset N_{\mbb R}$ be a complete 
fan with $\Si(1)\subseteq\Xi$. Let us denote with 
$p:U(\Si)\rar X_\Si$ and 
$r=pq:U(\Si)\times W_0^\times\rar X_\Si$ the projections 
and consider $\xi\in\Si(1)$. Then:

{\rm (i)} $\bigl(p_*{\cal O}_{U(\Si)}\bigr)^{c_{\Si}(f_\xi)}$ 
and $\bigl(r_*{\cal O}_{U(\Si)\times W_0^\times}
\bigr)^{c_{\Xi}(f_\xi)}$ are 
both isomorphic with $\cal O(D_\xi)$. If $\Si$ is regular, 
$\cal O(D_\xi)\cong U(\Si)\times_{c_\Si(f_\xi)}{\mbb C}
\cong \bigl(U(\Si)\times W_0^\times\bigr)
\times_{c_\Xi(f_\xi)}{\mbb C}$. The sheaves 
$\bigl(r_*{\cal O}_{U(\Si)\times W_0^\times}
\bigr)^{c_{\Xi}(f_j)}$, $j\in J$, are all isomorphic to 
$\cal O_{X_\Si}$.

{\rm (ii)} $p^*\bigl(p_*{\cal O}_{U(\Si)}\bigr)\rar U(\Si)$ 
is isomorphic with $\cal O_{U(\Si)}$, for all $\xi\in\Si(1)$. 
The $G_\Si$-action on $U(\Si)$ admits a canonical linearization 
in $p^*\cal O(D_\xi)$, such that for any open subset 
$O\subset X_\Si$ the sections of $\cal O(D_\xi)\rar O$ 
correspond to the $G_\Si$-invariant sections in 
$p^*\cal O(D_\xi)\rar p^{-1}O$. Under the isomorphism above, 
this linearization correponds to the action with character 
$c_\Si(f_\xi)$ of $G_\Si$ on $\cal O_{U(\Si)}$.

Similarly, $r^*\bigl(r_*{\cal O}_{U(\Si)\times W_0^\times}
\bigr)^{c_{\Xi}(f_\xi)}\bigr)$ is isomorphic with  
${\cal O}_{U(\Si)\times W_0^\times}$, for all $\xi\in\Xi$. 
The linearization of the $G_\Xi$-action corresponds, via 
this isomorphism to the action with weight $c_\Xi(f_\xi)$.
\end{proposition}

\begin{proof}
(i) We determine the sections of 
$\bigl(p_*{\cal O}_{U(\Si)}\bigr)^{c_{\Si}(f_\xi)}$ and 
$\bigl(
r_*{\cal O}_{U(\Si)\times W_0^\times}
\bigr)^{c_{\Xi}(f_\xi)}$ 
over the affine pieces $X_\si\subset X_\Si$. 
For $\si\in\Si$ a cone, we define the sets 
$$
M_\xi(\si)=
\begin{cases}
\si^\vee\cap M
=\{m\in M\mid\langle m,\xi'\rangle\geq 0,\;\forall \xi'\in\si(1)\}
\quad\text{if }\xi\not\in\si(1);
\\[1.5ex] 
\left\{
m\in M\;\biggl|
{\begin{array}{l}
\langle m,\xi'\rangle\geq 0,\;\forall \xi'\in\si(1)\sm\{\xi\}\\ 
\langle m,\xi\rangle\geq -1
\end{array}}\biggr.
\right\}
\quad\text{if }\xi\in\si(1).
\end{cases}
$$
We can give now explicitely describe the space of sections over 
$X_\si$ of the various sheaves which appear in the context:
$$
\begin{array}{l}
\Gamma
\bigl( 
X_\si, 
\bigl(p_*{\cal O}_{U(\Si)}\bigr)^{c_{\Si}(f_\xi)}
\bigr)
=
\bigoplus_{m\in M_\xi(\si)}
{\mbb C}\ {\cdot} Z_\xi\prod_{\xi'\in\Si(1)} 
Z_{\xi'}^{\langle m,\xi'\rangle},
\\[1.5ex]
\Gamma
\bigl( 
X_\si, 
\bigl(r_*{\cal O}_{U(\Si)\times W_0^\times}\bigr)^{c_{\Xi}(f_\xi)}
\bigr)
=
\bigoplus_{m\in M_\xi(\si)}
{\mbb C}\ {\cdot} Z_\xi\prod_{\rho\in\Xi} Z_\rho^{\langle m,\rho\rangle},
\\[1.5ex]
\Gamma
\bigl( 
X_\si, 
\bigl(r_*{\cal O}_{U(\Si)\times W_0^\times}\bigr)^{c_{\Xi}(f_j)}
\bigr)
=
\bigoplus_{m\in \si^\vee\cap M}
{\mbb C}\ {\cdot} Z_j\prod_{\rho\in\Xi} Z_\rho^{\langle m,\rho\rangle},
\\[1.5ex] 
\Gamma
\bigl( 
X_\si,{\cal O}_{X_\Si}(D_\xi)
\bigr)
=\bigoplus_{m\in M_\xi(\si)}
{\mbb C}\ {\cdot}\prod_{\xi'\in\Si(1)} 
Z_{\xi'}^{\langle m,\xi'\rangle}.
\end{array}
$$
The first three equalities are immediate consequences of 
the diagram \eqref{cd1}, while the fourth one can be found 
in \cite[page 66]{fu}. We claim that the homomorphism 
$$
\bigl(r_*{\cal O}_{U(\Si)}\bigr)^{c_{\Xi}(f_\eta)} 
\lar 
{\cal O}_{X_\Si}(D_\xi),\quad 
Z_\eta\prod_{\rho\in\Xi} Z_\rho^{\langle m,\rho\rangle}
\lmt 
\prod_{\xi'\in\Si(1)} Z_{\xi'}^{\langle m,\xi'\rangle} 
$$
defines a sheaf isomorphism for every 
$\eta\in\Xi=\Si(1)\cup J$, since every $m\in M$ is uniquely 
defined by its evaluation on the elements in $\Si(1)$.

For $j\in J$, the function 
$$
U(\Si)\times W_0^\times\lar {\mbb C},\quad 
(u,w)\lmt w_j
$$
is $G_\Xi$-equivariant with respect to the $c_\Xi(f_j)$-action 
on $\mbb C$, and therefore defines a global section 
$X_\Si\rar \bigl(p_*{\cal O}_{U(\Si)\times 
W_0^\times}\bigr)^{c_{\Xi}(f_j)}$ which vanishes nowhere. 
Moreover, the description of the section spaces given above 
above implies that $\bigl(p_*{\cal O}_{U(\Si)\times 
W_0^\times}\bigr)^{c_{\Xi}(f_j)}$ is locally free of 
rank one. It is therefore isomorphic to $\cal O_{X_\Si}$. 

(ii) We prove our statements for 
$p^*\bigl(p_*\cal O_{U(\Si)}\bigr)\rar U(\Si)$, 
since the other case is completely analogous. We distinguish 
again between two possibilities, according to whether $\xi$ 
belongs to $\si(1)$ or not. For $\xi\not\in\si(1)$,
$$
\begin{array}{ll}
{\mbb C}[U_\si]\otimes_{{\mbb C}[U_\si]^{G_\Si}}
\Gamma\bigl(
X_\si, \bigl(p_*{\cal O}_{U(\Si)}\bigr)^{c_{\Si}(f_\xi)}
\bigr)
&
=
{\mbb C}[U_\si]\otimes_{{\mbb C}[U_\si]^{G_\Si}}
{\mbb C}[U_\si]^{G_\Si}\cdot Z_\xi
\\[1.5ex]
&
={\mbb C}[U_\si]\otimes_{\mbb C}{\mbb C} Z_\xi,
\end{array}
$$
and we notice that in this case $Z_\xi$ is invertible 
in ${\mbb C}[U_\si]$. 

On the other hand, the smallest invariant open set 
$U_\si\subset U(\Si)$ with the property that $\xi\in\si(1)$ 
is $U_\xi:=U_{\mbb R_{\geq 0}\xi}$ (let us recall that 
${\mbb R}_{\geq 0}\xi$ is itself a cone of $\Si$). Then, 
if we fix some $m_0\in M_\xi({\mbb R}_{\geq 0}\xi)$ such 
that $\langle m_0,\xi\rangle=-1$, and write $Z^{m_0}:=
\prod_{\xi'\in\Si(1)}Z_{\xi'}^{\langle m_0.\xi'\rangle}$, 
we find that
$$
\begin{array}{ll}
{\mbb C}[U_{\xi}]\otimes_{{\mbb C}[U_\xi]^{G_\Si}}
\Gamma\bigl(
X_\xi, \bigl(p_*{\cal O}_{U(\Si)}\bigr)^{c_{\Si}(f_\xi)}
\bigr)
&
={\mbb C}[U_{\xi}]\otimes_{{\mbb C}[U_{\xi}]^{G_\Si}}
{\mbb C}[U_{\xi}]^{G_\Si}{\cdot}Z_\xi Z^{m_0}
\\[1.5ex] 
&
={\mbb C}[U_{\xi}]\otimes_{\mbb C}{\mbb C}Z_\xi Z^{m_0}.
\end{array}
$$
Again, since $\langle m_0,\xi\rangle=-1$, $Z_\xi Z^{m_0}$ is 
invertible in ${\mbb C}[U_{\xi}]$. Consequently, 
$$
p^*\bigl[\bigl(p_*\cal O_{U(\Si)}\bigr)^{c_{\Si}(f_\xi)}\bigr]
\xrar{\;\rm product\;}
{\cal O}_{U(\Si)}
$$
is an isomorphism of modules. Since $G_\Si$ acts on both $Z_\xi$ 
and $Z_\xi Z^{m_0}$ with character $c_\Si(f_\xi)$, we deduce that 
the induced linearization is given by the weight $c_\Si(f_\xi)$. 
\end{proof}

Let us remark that for any $\beta\in A_\Si$ and $\xi\in\Si(1)$ 
such that $c_\Si(f_\xi)=\beta$, 
$$
\Gamma
\bigl(
X_\Si,\bigl(p_*\cal O_{U(\Si)}\bigr)^{c_{\Si}(f_\xi)}
\bigr)
=R_\beta,
$$
the summand of degree $\beta$ of the homogeneous ring $R$. 
Moreover, $R_\beta$ decomposes into the direct sum 
$R_\beta'\oplus R_\beta''$, where $R_\beta'
=\oplus_{\{\xi\mid f_\xi\in c_\Si^{-1}(\beta)\}}\mbb CZ_\xi$ and 
$R_\beta''$ consists of sums of monomials of degree $\beta$ which 
are products of two or more variables. In the terminology of 
\cite[section 4]{cox}, $R_\beta'$ corresponds to the semi-simple 
roots of $X_\Si$, while $R_\beta''$ corresponds to the unipotent 
roots. 

The point of this discussion, which we are going to use in the 
next section, is that $W_\beta^\vee$ is canonically isomorphic 
to $R_\beta'$, and therefore is a direct summand in 
$\Gamma\bigl(X_\Si,
\bigl(p_*\cal O_{U(\Si)}\bigr)^{c_{\Si}(f_\xi)}\bigr)$. 

Let us consider now the ring 
$\tld R:={\mbb C}[Z_\rho\mid\rho\in\Xi]$. The group $G_\Xi$ 
acts on it, and the corresponding weights define a graduation 
$$
\tld R=\bigoplus_{\alpha\in A_\Xi}\tld R_\alpha.
$$
Similarly, for any $\alpha\in A_\Xi$ and $\xi\in\Xi$ 
such that $c_\Xi(f_\xi)=\alpha$,
$$
\Gamma
\bigl(
X_\Si,\bigl(r_*\cal O_{U(\Si)\times W_0^\times}\bigr)^{c_{\Xi}(f_\xi)}
\bigr)
=\tld R_\alpha,
$$
and $\tld R_\alpha=\tld R_\alpha'\oplus\tld R_\alpha''$, with  
$\tld R_\alpha'$ canonically isomorphic to $V_\alpha^\vee$. 

Now we fix $\xi\in \Si(1)$ and we wish to understand the 
relationship between $\tld R_{c_\Xi(f_\xi)}$ and 
$R_{c_\Si(f_\xi)}$. Let us define 
$$
M_\xi:=
\biggl\{
m\in M\biggl| 
\begin{array}{l}
\langle m,\xi\rangle\geq -1,\\ 
\langle m,\xi'\rangle\geq 0,\;\forall\,\xi'\in\Si(1)\sm\{\xi\}
\end{array}
\biggr.
\biggr\}.
$$
Then 
$$
\begin{array}{l} 
\Gamma(X_\Si,\tcl_\xi)
=\oplus_{m\in M_\xi} {\mbb C}\,
Z_\xi\cdot\prod_{\rho\in \Xi}Z_\rho^{\langle m,\rho\rangle}
\\ 
\Gamma(X_\Si,\crl L_\xi)
=\oplus_{m\in M_\xi} {\mbb C}\,
Z_\xi\cdot\prod_{\xi'\in \Si(1)}Z_{\xi'}^{\langle m,\xi'\rangle},
\end{array}
$$
and we notice that 
$$
\tld R_{c_\Xi(f_\xi)}=
\Gamma(X_\Si,\tcl_\xi)\srel{f}{\lar} \Gamma(X_\Si,\crl L_\xi)
=R_{c_\Si(f_\xi)}
$$ 
$$
Z_\xi\cdot\prod_{\rho\in \Xi}Z_\rho^{\langle m,\rho\rangle}
\lmt 
Z_\xi\cdot\prod_{\xi'\in \Si}Z_{\xi'}^{\langle m,\xi'\rangle}
$$
is an isomorphism of modules, as it should be, since 
$\tcl_\xi\cong\cal O_{X_\Si}(D_\xi)\cong\crl L_\xi$. 

We want to determine $f^{-1}R_{c_\Si(f_\xi)}'$. For shorthand, 
we write $\beta=c_\Si(f_\xi)$. First of all, 
$R_\beta'=\oplus_{\{\xi'\mid f_{\xi'}\in c_\Si^{-1}(\beta)\}}
{\mbb C}Z_{\xi'}$, by definition. We partition the indices as 
follows
$$
\{\xi'\mid f_{\xi'}\in c_\Si^{-1}(\beta)\}
=\bigsqcup_{\alpha\in b^{-1}(\beta)}\kern-1ex
\{\xi'\mid c_\Xi(f_{\xi'})=\alpha\}.
$$
\begin{lemma}{\label{lm:k}}
For every $\alpha\in b^{-1}(\beta)$, there are 
integers $(k_j(\alpha))_{j\in J}$ depending on $\alpha$, 
such that for all $f_{\xi'}\in c_\Xi^{-1}(\alpha)$ holds: 
$f^{-1}Z_{\xi'}=Z_{\xi'}\prod_{j\in J}Z_j^{k_j(\alpha)}$.
\end{lemma}

\begin{proof} Let us fix such $\alpha$ and $\xi_0'$: then 
$c_\Si(f_{\xi_0'})=\beta=c_\Si(f_\xi)$, and therefore 
we find a unique $m_{\xi_0'}\in M$ such that 
$\langle m_{\xi_0'},\xi\rangle=-1$, 
$\langle m_{\xi_0'},\xi_0'\rangle=+1$, and all the other 
evaluations on vectors in $\Si(1)\sm \{\xi,\xi_0'\}$ vanish. 
Then $Z_{\xi_0'}=Z_\xi\cdot\prod_{\xi'\in\Si(1)}
Z_{\xi'}^{\langle m,\xi'\rangle}$ and 
$f^{-1}Z_{\xi_0'}=Z_\xi\cdot\prod_{\rho\in \Xi}
Z_\rho^{\langle m_{\xi_0'},\rho\rangle}
=Z_{\xi_0'}\prod_{j\in J}Z_j^{\langle m_{\xi_0'},j\rangle}$. 
We define $k_j(\alpha):=\langle m_{\xi_0'},j\rangle$, and we 
must check that these numbers do not depend on the choice of 
$\xi_0'$ satisfying $c_\Xi(f_{\xi_0'})=\alpha$. If $\xi_0''$ 
is another vector with the same property, we find a unique 
$m_{\xi_0'\xi_0''}\in M$ which evaluates $-1$ on $\xi_0'$, 
$+1$ on $\xi_0''$, and zero on the other vectors in 
$\Xi\sm\{\xi_0',\xi_0''\}$. From the unicity of the elements 
$\{m_{\xi'}\}_{\xi'}$, we deduce that 
$m_{\xi_0''}=m_{\xi_0'}+m_{\xi_0'\xi_0''}$, and the evaluations 
of it on $j\in J$ coincide with those of $m_{\xi_0'}$.
\end{proof}

The lemma can now be reformulated as
\begin{align}{\label{eq:f-1}}
f^{-1}R_\beta'=\kern-1ex
\bigoplus_{\alpha\in b^{-1}(\beta)}\kern-1ex 
\tld R_{\alpha}\otimes_{\mbb C}{\mbb C}
\mbox{$\prod_{j\in J}Z_j^{k_j(\alpha)}$}
=\kern-1ex
\bigoplus_{\alpha\in b^{-1}(\beta)}\kern-1ex 
V_{\alpha}^\vee\otimes_{\mbb C}{\mbb C}
\mbox{$\prod_{j\in J}{(V_j^\vee)}^{\otimes k_j(\alpha)}$}, 
\end{align}
and this is the form that we will need subsequently.


\section{The construction}{\label{sct:constr}}

This section is devoted presenting a method for constructing 
families of complete toric varieties, parameterized by some 
base space. The idea behind the construction is to replace 
${\mbb C}^\Xi$, which should be seen as the as the total space 
of a vector bundle over a point, with an arbitrary vector bundle 
$\euf V\rar S$ on which a torus acts. The result of the 
construction should be a locally trivial family of toric 
varieties, parameterized by the base $S$. 

Let us make things more precise: the main ingredients are 
again the lattice $N$ and the collection of vectors 
$\Xi\subset N$ having the properties (i)--(iii) enumerated 
in the previous section. Associated to this data we have 
the exact sequence 
$$
1\lar G_\Xi\lar {({\mbb C}^*)}^\Xi\lar T_N\lar 1.
$$
Corresponding to the induced $G_\Xi$-action, we obtain the 
decomposition of ${\mbb C}^\Xi$ into the $G_\Xi$-isotypic 
components
$$
{\mbb C}^\Xi=\bigoplus_{\alpha\in A_\Xi}V_\alpha,\quad 
\dim V_\alpha=r_\alpha.
$$
The second ingredient for the construction is a separated, 
noetherian scheme $S$ and a collection of arbitrary vector 
bundles (locally free sheaves) $\euf V_\alpha\rar S$ of rank 
$r_\alpha$ respectively. We define the vector bundle
\begin{align}{\label{xi-bdl}}
\euf V:=\bigoplus_{\alpha\in A_\Xi}
\unl{\mbb C}_\alpha\otimes \euf V_\alpha,
\end{align}
where $\unl{\mbb C}_\alpha$ stands for the trivial vector 
bundle over $S$ on which $G_\Xi$ acts by the character 
$\alpha$. For this choice of the vector bundle 
$\euf V\rar S$, we  obtain fibrewise the same $G_\Xi$-action 
as the one described in the preceeding section. 

At this point we recall a result concerning the structure 
of vector bundles with a linear group action, which will be 
needed later on.\smallskip 

\nit\unl{\sl Result}\;\cite[prop. 1]{kraft}\;
{\it Let $G$ be a reductive group and $\euf V\rar S$ a 
$G$-vector bundle, where $G$ acts trivially on $S$. 

\rm (i) \it $\euf V\rar S$ is locally trivial, as a $G$-vector 
bundle, in the Zariski topology. In particular, if $S$ is 
connected, all the representations $G\rar Gl(\euf V_s)$, 
$s\in S$, are equivalent.

\rm (ii) \it There is an isomorphism of $G$-vector bundles 
$$
\euf V\srel{\cong}{\lar}\bigoplus_\omega M_\omega\otimes\euf V_\omega,
$$
where the $M_\omega$'s are simple $G$-modules and the 
$\euf V_\omega$'s are vector bundles with trivial $G$-action.
The direct summands $M_\omega\otimes \euf V_\omega$ are called 
the $\omega$-isotypic components of $\euf V$, and they are 
uniquely defined.
}\smallskip

We choose now a complete fan $\Si\subset N_{\mbb R}$ with 
$\Si(1)\subseteq\Xi$ and we want to construct a locally 
trivial fibration over $S$ with all the fibres isomorphic 
to the toric variety $X_\Si$. Since $G_\Si$ is a subgroup 
of $G_\Xi$, the vector bundle $\euf V$ gets an induced 
$G_\Si$-action. For every $\alpha\in A_\Xi$, $G_\Si$ acts on 
$\unl{\mbb C}_\alpha\otimes\euf V_\alpha$ by the character 
$b(\alpha)\in A_\Si$, where the homomorphism $b:A_\Xi\rar 
A_\Si$ was defined in \eqref{cd1}. We deduce that the 
$G_\Si$-isotypical decomposition of $\euf V$ is
\begin{align}{\label{si-bdl}}
\euf V=\bigoplus_{\beta\in A_\Si}\euf W_\beta
=\biggl(
\bigoplus_{\beta\in A_\Si\sm\{0\}}\euf W_\beta
\biggr)\oplus \euf W_0
=:\euf W_\Si\oplus\euf W_0. 
\end{align}
In this formula, for each $\beta\in A_\Si$, we have denoted 
$\euf W_\beta:=\bigoplus_{\alpha\in b^{-1}(\beta)}
\unl{\mbb C}_\alpha\otimes\euf V_\alpha$.

\begin{lemma}{\label{lm:0}}
$\euf W_0=\bigoplus_{\alpha\in\ker(b)}\unl {\mbb C}_\alpha
\otimes \euf V_\alpha\rar S$ is a vector bundle whose rank 
equals $\#(\Xi\sm\Si(1))$. Moreover, each of the 
$\unl{\mbb C}_\alpha\otimes \euf V_\alpha$'s has rank equal 
one.
\end{lemma}

\begin{proof}
By the first part of the result above, $\euf V$ is locally 
trivial, as a $G_\Si$-bundle, in the Zariski topology. 
Therefore, locally over $S$, we are in the situation of 
the lemma \ref{lm:triv-subsp} and the conclusion follows.
\end{proof}

\begin{definition}{\label{def:prim}}
We shall write 
$\euf W_0^\times:=\bigoplus_{j\in\Ker(b)}
\bigl[\bigl(\unl{\mbb C}_{j}\otimes 
\euf V_j\bigr)\sm\{0\}\bigr]$. 
For a $\Si$-primitive family $\pi\subseteq\Si(1)$, we define 
the subvector bundle
$$
\euf W(\pi):=\Ker\biggl[
\pr_\pi:\euf W_\Si\lar\bigoplus_{\beta\subset\pi}
\euf W_\beta\biggr]
\subset\euf W_\Si.
$$ 
The meaning of the notation $\beta\subset\pi$ was defined in 
\ref{def:incl}.

We denote by $\euf Z_\Si$ the union of these bundles as $\pi$ 
varies among the $\Si$-primitive families
$$
\euf Z_\Si:=\bigcup_{\pi\text{ primitive}}\kern-1ex\euf W(\pi)
\subset\euf W_\Si.
$$
\end{definition}

\begin{theorem}{\label{the-constr}}
Let us consider $\Si\subset N_{\mbb R}$ a complete fan 
with $\Si(1)\subseteq \Xi$, and define the open subsets 
$\crl U(\Si):=\euf W_\Si\sm\euf Z_\Si\subset \euf W_\Si$ 
and 
$\tcu(\Si):=\crl U(\Si)\times_S\euf W_0^\times\subset\euf V$. 
Then the categorical quotients $\tcu_\Si/G_\Xi$ and 
$\crl U(\Si)/G_\Si$ does exist. They are locally trivial 
fibrations over $S$, with the fibres isomorphic to the 
toric variety $X_\Si$ defined by the fan $\Si$. When $\Si$ 
is moreover simplicial, the above quotients are geometric.
\end{theorem}

\begin{proof}
We shall prove the statement concerning 
$\tcu(\Si)$, 
the proof for ${\crl U}(\Si)$ being completely identical. 

The first thing to notice is that $\tcu(\Si)$ is 
$G_\Xi$-invariant, since ${\crl U}(\Si)$ is $G_\Si$-invariant. 
We know that the $G_\Xi$-bundle $\euf V$ is locally trivial 
in the Zariski topology. In order to construct the desired 
quotient we can restrict ourselves to open subsets of $S$: 
let us choose some open subset $O\subset S$ over which 
$\euf V$ is trivializable. It follows from \eqref{eq:equiv} that 
$\tcu(\Si)|_O\cong O\times (U(\Si)\times W_0^\times)$, 
where $U(\Si)={\mbb C}^{\Si(1)}\sm Z_\Si$ (see section 
\ref{sct:rmk}), the isomorphism being $G_\Xi$-equivariant. We use 
now proposition \ref{prop:inv-quot} and deduce the existence of 
the categorical quotient $\tcu(\Si)|_O/G_\Xi$, which 
is isomorphic with $O\times X_\Si$. From the universality 
property of the categorical quotient, it follows that, as 
$O\subset S$ varies, we can paste these quotients together. 
Therefore the categorical quotient $\tcu(\Si)/G_\Si$ 
exists, and it is a $X_\Si$-fibration over the base $S$.

When $\Si$ is simpicial, we know that 
$(U(\Si)\times W_0^\times)/G_\Xi$ is a geometric quotient, 
and consequently, $\tcu(\Si)/G_\Xi$ is still geometric.
\end{proof}

\begin{remark}
The most classical example of such a family is the 
projectivized bundle of an arbitrary vector bundle 
over some base variety. 

We must clearly point out that the two quotients 
$\tcu(\Si)/G_\Xi$ and $\crl U(\Si)/G_\Si$ are 
{\em not} isomorphic in general. Of course, the case 
$\Si(1)=\Xi$ is an exception, since then the two families 
coincide. The reason for obtaining in general different 
fibrations is due to the fact that the isomorphism between 
$\bigl(U(\Si)\times W_0^\times\bigr)/G_\Xi$ and $U(\Si)/G_\Si$ 
described in the proposition \ref{prop:inv-quot} was not natural. 
When $\Si$ is simplicial, we will see in the next section that 
we can distinguish them even at the cohomological level. 
\end{remark}

The families $\tcx_\Si$ and $\crl X_\Si$ are related 
in the following fashion: let us take the quotient 
$\tcu(\Si)/G_\Xi$ in two steps, first for the 
$G_\Si$-action and then for the remaining $G_J=G_\Xi/G_\Si$-action. 
We deduce that 
$$
\tcx_\Si
=\bigl(\crl X_\Si\times_S \euf W_0^\times\bigr)\bigl/G_J\bigr., 
$$
and this quotient is always geometric since $\euf W_0^\times\rar S$ 
is a principal $G_J$-bundle. 

We shall denote with $\tcx_\Si\rar S$ and 
$\crl X_\Si\rar S$ respectively the two families constructed 
above. For a given pair $(N,\Xi)$ and a given collection of 
vector bundles $\{\euf V_\alpha\rar S\}_\alpha$, we remark that 
both $\tcx_\Si$ and $\crl X_\Si$ depend on the choice of the 
complete fan $\Si\subset N_{\mbb R}$ with $\Si(1)\subseteq \Xi$. 
There are several ways of choosing complete fans satisfying this 
property, and the detailed account about this issue can be found 
in \cite{op}. The moral of the discussion is that one may construct 
several different fibrations out of a given vector bundle 
$\euf V\rar S$ with a given $G_\Xi$-action. 

We are going now to take advantage of the invariant, 
coordinate-free description of the sheaves 
$\cal O_{X_\Si}\rar X_\Si$, and define their analogues 
for families of toric varieties. Before doing this, we recall 
from lemma \ref{lm:triv-subsp} that for $j\in J$, $\euf V_j\rar S$ 
is a line bundle.

\begin{proposition}{\label{prop:lxi}}
Let $\Si$ be a complete fan with $\Si(1)\subseteq\Xi$. 
For each $\xi\in\Si(1)$ we define 
$$
\tcl_\xi:=
{\bigl(
r_*\cal O_{\tcu(\Si)}
\bigr)}^{c_\Xi(f_\xi)}\lar\tcx_\Si
\text{ and }
\crl L_\xi:=
{\bigl(
p_*\cal O_{\crl U(\Si)}
\bigr)}^{c_\Si(f_\xi)}\lar\crl X_\Si,
$$
where $r:\tcu(\Si)\rar\tcx_\Si$ and 
$p:\crl U(\Si)\rar\crl X_\Si$ are the projections. 

\nit {\rm (i)}  $\tcl_\xi$ and $\crl L_\xi$ are coherent, 
torsion free sheaves of rank one. Their restiction to 
the fibres of the corresponding families are isomorphic 
with $\cal O_{X_\Si}(D_\xi)$. Moreover, $\tcl_\xi$ 
and $\crl L_\xi$ are locally free outside a singular set of 
codimension at least two.

When $\Si$ is regular, both sheaves are locally free:
$$
\tcl_\xi=\tcu(\Si)\times_{c_\Xi(f_\xi)}\mbb C 
\text{ and }
{\crl L}_\xi={\crl U}(\Si)\times_{c_\Si(f_\xi)}\mbb C. 
$$

\nit {\rm (ii)} $r^*\tcl_\xi\rar \tcu(\Si)$ and 
$p^*\crl L_\xi\rar{\crl U}(\Si)$ are respectively isomorphic 
to $\cal O_{\tcu(\Si)}$ and $\cal O_{{\crl U}(\Si)}$, 
and they carry natural $G_\Xi$ and $G_\Si$-linearizations. 
These linearizations correpond respectively, under 
the isomorphisms, to the actions with weights $c_\Xi(f_\xi)$ and 
$c_\Si(f_\xi)$.

\nit {\rm (iii)} For $j\in J$, $\tcl_j:={\bigl(r_*\cal O_{\tcu(\Si)}
\bigr)}^{c_\Xi(f_j)}\lar\tcx_\Si$ is isomorphic to 
$\tld\phi^*\euf V_j^\vee$, where $\tld\phi:\tcx_\Si\rar S$ is the 
projection. Moreover, $r^*\tcl_j\cong \cal O_{\tcx_\Si}$, and 
the induced linearization of $G_\Xi$ correponds to the action with 
weight $c_\Xi(f_j)$. 
\end{proposition}

\begin{proof}
(i) Again, we shall treat only the case of $\tcl_\xi$. 
Let us consider an open subset $O\subset S$ having the property 
that $\euf V|_O\rar O$ is trivializable. The proposition 
\ref{prop:coord-f} implies that over such trivializing open 
subsets
$$
\tcl_\xi|_O
\cong 
O\times\cal O_{X_\Si}(D_\xi)
\lar 
O\times{X_\Si}.
$$
Is is clearly a torsion-free sheaf, and its rank equals one. 
Since $\cal O_{X_\Si}(D_\xi)\rar X_\Si$ is coherent, and the 
coherence property is local in nature, we deduce the coherence 
of $\tcl_\xi\rar \tcx_\Si$. Moreover, since 
$X_\Si$ is a normal variety (according to \cite[page 29]{fu}), 
we deduce that $\cal O_{X_\Si}(D_\xi)\rar X_\Si$ is locally 
free outside a closed subset of codimension at least two in 
$X_\Si$. The isomorphism above implies that the same holds for 
$\tcl_\xi\rar\tcx_\Si$.

(ii) The germ of $r^*\tcl_\xi$ at a point $\tld u\in \tcu(\Si)$ 
is
$$
\bigl(\tcl_\xi\bigr)_{\tld u}=
\lim_{r(\tld u)\in O}
\biggl(
\cal O_{\tcu(\Si)}(r^{-1}O)
\otimes_{\cal O_{\kern.5ex\tld{\kern-.5ex\crl X}_\Si}(O)}
\bigl(\cal O_{\tcu(\Si)}(r^{-1}O)\bigr)^{c_\Xi(f_\xi)}
\biggr),
$$
and we find the natural homomorphism 
$$
\bigl(\tcl_\xi\bigr)_{\tld u}
\xrar{\ \rm product\ } 
\lim_{r(\tld u)\in O}\cal O_{\tcu(\Si)}(r^{-1}O)
\lar 
\cal O_{\tcu(\Si),\tld u}.
$$
In other words, there is a well-defined a homomorphism 
$r^*\tcl_\xi\rar \cal O_{\tcu(\Si)}$, and we have checked 
in proposition \ref{prop:coord-f} that it is actually an 
isomorphism. The statement about the linearizations is 
proved in the same place.

(iii) For the last part, let us choose a trivializing covering 
$\{O_\nu\}_\nu$ of $S$, for $\euf V_j$, having the property that 
there are non-vanishing sections 
$s_\nu:\cal O_S(O_\nu)\rar \euf V_j|_{O_\nu}$. Then $\euf V_j$ 
is associated to the co-cycle $\{s_\mu^{-1}s_\nu\}_{\nu,\mu}$. 
Let $\pr_s:\tcu(\Si)\rar S$ be the projection. Then we have a 
nowhere vanishing section
$$
s_j:\tcu(\Si)\lar \pr_S^*\euf V_j,\quad (u,w)\lmt 
\bigl((u,w),w_j\bigr).
$$
For each $\nu$, the composition
$$
s_{j,\nu}:\tcu(\Si)|_{O_\nu}
\xrar{\;s_j\;}
\pr_S^*\euf V_j|_{O_\nu}
\xrar{\;\pr_S^*s_\nu^{-1}\;}
\pr_S^*\cal O_S(O_\nu)=\cal O_{\tcu(\Si)}|_{O_\nu}
$$
is $G_\Xi$-equivariant and non-vanishing, and therefore defines 
a trivializing section $\tcx_\Si|_{O_\nu}\rar \tcl_j|_{O_\nu}$. 
We deduce that the co-cycle defining $\tcl_j$ is 
${\{s_{j,\mu}^{-1}s_{j,\nu}\}}_{\nu,\mu}$, with 
$$
s_{j,\mu}^{-1}s_{j,\nu}=\tld\phi^*(s_\mu^{-1}s_\nu)
=\tld\phi^*\bigl((s_\mu^{-1}s_\nu)^{-1}\bigr),
$$
so that $\tcl_j\cong\tld\phi^*\euf V_j^{-1}$. The statement 
about the linearizations follows from the same proposition 
\ref{prop:coord-f}.
\end{proof}

It is a well known fact that for a projective bundle 
${\mbb P}(\euf W)\rar S$, the direct image on $S$ of 
$\cal O_{{\mbb P}(\euf W)}(1)$ is isomorphic to $\euf W^\vee$. 
We wish to investigate now what happens in our more general 
setting, that is we wish to compute the direct images on $S$ 
of the sheaves $\tcl_\xi$ and $\crl L_\xi$. 

Let us denote by $\euf P_\Xi$ and $\euf P_\Si$ respectively 
the principal bundles of frames in $\euf V\rar S$ and 
$\euf W_\Si\rar S$: the structure group of $\euf P_\Xi$ is 
$\cal G(\euf P_\Xi):=\prod_{\alpha\in A_\Xi}Gl(V_\alpha)$, 
while the structure group of $\euf P_\Si$ is 
$\cal G(\euf P_\Si):=\prod_{\beta\in A_\Si}Gl(W_\beta)$. 

\begin{proposition}{\label{prop:d-image}}
{\rm (i)}  Consider $\xi\in\Si(1)$ and denote 
$c_\Xi(f_\xi)=\alpha$ and $c_\Si(f_\xi)=\beta$. Then 
$$
\tld\phi_*\tcl_\xi\cong 
\euf P\times_{\cal G(\euf P_\Xi)}\tld R_\alpha,
\quad\text{for }\tld\phi:\tcx_\Si\rar S,
$$
and contains 
$$
\bigoplus_{\alpha'\in b^{-1}(\beta)}
\biggl[
\euf V_{\alpha'}^\vee
\otimes
\bigl(
\otimes_{j\in J}\euf V_j^{-k_j(\alpha')}
\bigr)
\biggr]
$$ 
as a direct summand. For $j\in J$, 
$\tld\phi_*\tcl_j\cong \euf V_j^\vee$.

\nit {\rm (ii)} For every $\beta\in A_\Si$ and $\xi\in \Si(1)$ 
with $c_\Si(f_\xi)=\beta$,
$$
\phi_*\crl L_\xi\cong 
\euf P_\Si\times_{\cal G(\euf P_\Si)}\tld R_\alpha,
\quad\text{for }\phi:\crl X_\Si\rar S,
$$
and contains $\euf W_\beta^\vee$ as a direct summand. 

\nit In particular, if $X_\Si$ has no unipotent roots, 
$\tld\phi_*\tcl_\xi$ and $\phi_*\crl L_\xi$ are isomorphic  
respectively to $\bigoplus_{\alpha'\in b^{-1}(\beta)}
\bigl[
\euf V_{\alpha'}^\vee
\otimes
\bigl(
\otimes_{j\in J}\euf V_j^{-k_j(\alpha')}
\bigr)
\bigr]$ and $W_\beta^\vee$.
\end{proposition}

\begin{proof} Since for all $j\in J$, 
$\tcl _j\cong\tld\phi^*\euf V_j^\vee$ and 
$X_\Si$ is connected, $\phi_*\tcl_j\cong \euf V_j^\vee$.

Let us choose a covering $\{O_\nu\}_\nu$ of $S$ which 
is simultaneously trivializing for the bundles 
$\{\euf V_\alpha\}_{\alpha\in A_\Xi}$. It is also a trivializing 
covering for the bundles $\euf P_\Xi$ and $\euf P_\Si$, and the 
corresponding transition maps $\{\tld g_{\mu\nu}\}_{\mu,\nu}$ 
and $\{g_{\mu\nu}\}_{\mu,\nu}$ are elements in 
$\cal G(\euf P_\Xi)$ and $\cal G(\euf P_\Si)$ respectively. 
We observe that these transition functions define automorphisms 
of $\tld U(\Si):=U(\Si)\times W_0^\times$ and $U(\Si)$ 
respectively, which in turn induce automorphisms of $X_\Si$ 
({\it cf.} \cite[section 4]{cox}). We notice further that since 
$\{\tld g_{\mu\nu}\}$ and $\{g_{\mu\nu}\}$ respect the graduation 
of the rings $\tld R$ and $R$ respectively, they induce automorphisms 
of $\bigl(r_*\cal O_{\tld U(\Si)}\bigr)^{c_\Xi(f_\xi)}$ 
and $\bigl(p_*\cal O_{U(\Si)}\bigr)^{c_\Si(f_\xi)}$ such that the 
diagrams 
$$
\xymatrix@C=4ex{
\bigl(r_*\cal O_{\tld U(\Si)}\bigr)^{c_\Xi(f_\xi)}
\ar[r]^-{\tld g_{\mu\nu}}\ar[d]
&
\bigl(r_*\cal O_{\tld U(\Si)}\bigr)^{c_\Xi(f_\xi)}
\ar[d]
\ar@{}[dr]|{\hbox{and}}
&
\bigl(p_*\cal O_{U(\Si)}\bigr)^{c_\Si(f_\xi)}
\ar[r]^-{g_{\mu\nu}}\ar[d]
&
\bigl(p_*\cal O_{U(\Si)}\bigr)^{c_\Si(f_\xi)}
\ar[d]
\\ 
X_\Si\ar[r]^-{\tld g_{\mu\nu}}
&
X_\Si
&
X_\Si\ar[r]^-{g_{\mu\nu}}
&
X_\Si
}
$$
commute. It follows then that 
$\tcx_\Si=\euf P_\Xi\times_{\cal G(\euf P_\Xi)}X_\Si$ and 
$\crl X_\Si=\euf P_\Si\times_{\cal G(\euf P_\Si)}X_\Si$. 
Moreover, above the trivializing open sets we have the 
isomorphisms 
$\tld\phi_*\tcl_\xi|_{O_\nu}\cong O_\nu\times \tld R_\alpha$ 
and $\phi_*\crl L_\xi|_{O_\nu}\cong O_\nu\times R_\beta$, 
and from what we said we deduce that $\tld\phi_*\tcl_\xi$ 
and $\phi_*\crl L_\xi$ are the the vector bundles associated 
to $\euf P_\Xi$ and $\euf P_\Si$ for the representations 
$\cal G(\euf P_\Xi)\rar Gl(\tld R_\alpha)$ and 
$\cal G(\euf P_\Si)\rar Gl(R_\beta)$ respectively. This 
proves the first half of the statements (i) and (ii). 

For the second half, we have proved in the end of section 
\ref{sct:rmk} that $W_\beta^\vee=R_\beta'$ is direct summand 
in $R_\beta$, and the equation \eqref{eq:f-1} shows that 
$f^{-1}R_\beta'=
\bigoplus_{\alpha\in b^{-1}(\beta)}
V_{\alpha}^\vee\otimes_{\mbb C}{\mbb C}
\mbox{$\prod_{j\in J}{(V_j^\vee)}^{\otimes k_j(\alpha)}$}$ is 
direct summand in $\tld R_\alpha$. Moreover, these sub-vector 
spaces are stable under the actions of the structure groups 
$\cal G(\euf P_\Si)$ and $\cal G(\euf P)$ respectively. This 
proves the rest of the proposition.
\end{proof}

Let us assume now that $S$ is projective. We want to prove 
that in this case the sheaves $\tcl_\xi\rar \tcx_\Si$ and 
$\crl L_\xi\rar \crl X_\Si$ define elements $[\tcl_\xi]
\in A^1(\tcx_\Si)$ and $[\crl L_\xi]\in A^1(\crl X_\Si)$. 
We have proved in proposition \ref{prop:lxi} that both 
$\tcl_\xi$ and $\crl L_\xi$ are locally free outside 
$\tcx_{\Si,\rm bad}\hra \tcx_\Si$ and 
$\crl X_{\Si,\rm bad}\hra\crl X_\Si$ respectively, 
whose codimension is at least two. Since $S$ is projective 
we can choose a sufficiently ample line bundle $A\rar S$ 
having the property that $\tld\phi_*\tcl_\xi\otimes A$, 
$\phi_*\crl L_\xi\otimes A$ and $A$ are all globally 
generated over $S$. Let us choose regular sections 
$$
\begin{array}{l}
\tld s\in\Gamma(S,\tld\phi_*\tcl_\xi\otimes A)=
\Gamma(\tcx_\Si,\tcl_\xi\otimes \tld\phi^*A),
\\ 
s\in\Gamma(S,\phi_*\crl L_\xi\otimes A)=
\Gamma(\crl X_\Si,\crl L_\xi\otimes \tld\phi^*A)\text{ and } 
s_A\in\Gamma(S,A), 
\end{array}
$$
which do not identically vanish on any irreducible 
component of $S$. 

Using them we can associate to the restrictions 
$\tcl_\xi|_{\tcx_\Si\sm \tcx_{\Si,\rm bad}}$ and 
$\crl L_\xi|_{\crl X_\Si\sm \crl X_{\Si,\rm bad}}$ the 
classes $[\tcl_\xi|_{\tcx_\Si\sm \tcx_{\Si,\rm bad}}]
:=[\text{Zero}(\tld s)]-[\text{Zero}(\tld\phi^*s_A)]\in 
A^1(\tcx_\Si\sm \tcx_{\Si,\rm bad})$ and 
$[\crl L_\xi|_{\crl X_\Si\sm \crl X_{\Si,\rm bad}}] 
:=[\text{Zero}(s)]-[\text{Zero}(\phi^*s_A)]\in 
A^1(\crl X_\Si\sm \crl X_{\Si,\rm bad})$. These classes 
do not depend on the choice of the sections $\tld s, s$ 
and $s_A$. We recall now that 
$\tcx_{\Si,\rm bad}\hra \tcx_\Si$ 
and $\crl X_{\Si,\rm bad}\hra \crl X_\Si$ both have 
codimension at least two, and therefore 
$$
A^1(\tcx _\Si)\lar A^1(\tcx_\Si\sm \tcx_{\Si,\rm bad})
\text{ and }
A^1(\crl X_\Si)\lar A^1(\crl X_\Si\sm \crl X_{\Si,\rm bad})
$$
are isomorphisms ({\it cf.} \cite[proposition 1.8]{ful}). 
We should further notice that these classes depend only 
on $c_\Xi(f_\xi)\in A_\Xi$ and $c_\Si(f_\xi)\in A_\Si$, 
and not of $\xi$ itself.


\section{Cohomology}{\label{sct:coh}}

In this section we determine the cohomology rings of the 
families of toric varieties constructed in the previous 
section, in the case when the fibre is simplicial and 
complete. This means that our context is as follows: $N$ 
is a lattice, $\Xi\subset N$ is a finite set satisfying 
the properties (i)--(iii), and $\Si\subset N_{\mbb R}$ is 
a complete and simplicial fan with $\Si(1)\subseteq\Xi$. 
Associated to this data we have the commutative diagrams 
\eqref{cd1} and \eqref{cd2}.

The goal of this section is to determine the cohomology 
rings of the families $\tcx_\Si$ and ${\crl X}_\Si$ constructed 
in the previous section. In particular, we shall see that 
they are usually not isomorphic, and can be distinguished 
by cohomological means.

We start with some well-known facts concerning classifying 
spaces for groups. 

\begin{definition}{\label{cl-space}}
We denote by $E$ a contractible space on which the group 
$\cxi$ acts freely. It is uniquely defined up to homotopy 
equivalence. The quotients 
$$
\BG_\Xi:=E/G_\Xi\text{ and } \BG_\Si:=E/G_\Si
$$
are called classifying spaces for $G_\Xi$ and $G_\Si$ 
respectively.
\end{definition}

\begin{lemma}{\label{lm:bg}}
$$
\begin{array}{ll}
H^*(\BG_\Xi;{\mbb Q})
&
=\Sym^\bullet(A_\Xi)_{\mbb Q}
\cong{\mbb Q}[f_\xi\mid\xi\in\Xi]/\LR_\Xi;
\\[1.5ex] 
H^*(\BG_\Si;{\mbb Q})
&
=\Sym^\bullet(A_\Si)_{\mbb Q}
\cong{\mbb Q}[f_\xi\mid\xi\in\Si(1)]/\LR_\Si
\\[1.5ex] 
&
\cong H^*(\BG_\Xi;{\mbb Q})\bigl/
\bigl(\langle f_j\mid j\in J\rangle+\LR_\Xi\bigr)
\bigr.,
\end{array}
$$
where 
$\LR_\Xi:=\langle a_\Xi(m)\mid m\in M\rangle\subset 
{\mbb Q}[f_\xi\mid\xi\in\Xi]$, and similarly for 
$\LR_\Si$. We call them the ideals of linear relations. 
Moreover, if $A_\Xi$ has no torsion, both isomorphisms 
hold with integral coefficients.
\end{lemma}

For shorthand, we write $\hg_{G_\Xi}$ and $\hg_{G_\Si}$ 
respectively for the rings $H^*(\BG_\Xi;-)$ and 
$H^*(\BG_\Si;-)$, where the coefficient ring is $\mbb Z$ 
or $\mbb Q$, depending on whether $A_\Xi$ has torsion or has 
not.

Let us remark at this point that if there exists a regular 
fan $\Si$ with $\Si(1)\subset\Xi$, then $A_\Xi$ has no torsion 
automatically. Indeed, $b:(A_\Xi)_{\rm tor}\rar 
(A_\Si)_{\rm tor}$ is always an isomorphism, and for regular 
fans $(A_\Si)_{\rm tor}=0$.

\begin{proof}
The torsions $(A_\Xi)_{\rm tor}\subset A_\Xi$ and 
$(A_\Si)_{\rm tor}\subset A_\Si$ are finite abelian groups and 
$b:(A_\Xi)_{\rm tor}\rar (A_\Si)_{\rm tor}$ is an isomorphism. 
The quotients $A_\Xi^o:=A_\Xi/(A_\Xi)_{\rm tor}$ and 
$A_\Si^o:=A_\Si/(A_\Si)_{\rm tor}$ are free $\mbb Z$-module and 
their duals $G^o_\Xi$ and $G^o_\Si$ are respectively the connected 
components of the identity in $G_\Xi$ and $G_\Si$. These latter 
groups are respectively isomorphic to the direct products 
$G^o_\Xi\times F$ and $G^o_\Si\times F$, with 
$F:=\Hom_{\mbb Z}((A_\Si)_{\rm tor},{\mbb C}^*)$. The 
classifying spaces for $G^o_\Xi$ and $G^o_\Xi$ are obtained as 
$\BG^o_\Xi:=E/G^o_\Xi$ and $\BG^o_\Si:=E/G^o_\Si$ respectively, 
and they fit in the commutative diagram
$$
\xymatrix{
E\bigl/\cxi\bigr.={\rm B}\cxi 
&
\ar[l] \BG^o_\Xi=E/G^o_\Xi
\\
\ar[u] 
E\bigl/\csi\bigr.={\rm B}\csi 
&
\ar[l] \ar[u]\BG^o_\Si=E/G^o_\Si.
}
$$
It is known that for a complex torus, the cohomology of its 
classifying space is naturally isomorphic to the symmetric 
algebra of its character group. In our case, the previous 
diagram induces
\begin{align}{\label{cd3}}
\xymatrix{
H^*_\cxi
=\Sym^\bullet\bigl({\mbb Z}^\Xi\bigr)
\ar@{->>}[r]^-{u_\Xi}\ar@{->>}[d]
&
H^*_{G_\Xi^o}
=\Sym^\bullet\bigl({\cal X}^\star(G^o_\Xi)\bigr)
=\Sym^\bullet\bigl(A_\Xi^o\bigr)
\ar@{->>}[d]^-b
\\
H^*_\csi
=\Sym^\bullet\bigl({\mbb Z}^{\Si(1)}\bigr)
\ar@{->>}[r]^-{u_\Si}
&
H^*_{G_\Si^o}
=\Sym^\bullet\bigl({\cal X}^\star(G^o_\Si)\bigr)
=\Sym^\bullet\bigl(A_\Si^o\bigr).
}
\end{align}
The homomorphisms $u_\Xi$ and $u_\Si$ are surjective and 
$$
\Ker(u_\Xi)=\langle\Ker(c^o_\Xi)\rangle,\quad 
\Ker(u_\Si)=\langle\Ker(c^o_\Si)\rangle,
$$
with $c^o_\Xi:{\mbb Z}^\Xi\lar A_\Xi^o$ and similarly 
$c^o_\Si$. If $(A_\Xi)_{\rm tor}=0$, we are done because 
in this case $G_\Xi$ and $G_\Si$ coincide with $G^o_\Xi$ 
and $G^o_\Si$ respectively, and 
$\langle\Ker(c^o_\Xi)\rangle=\LR_\Xi$, 
$\langle\Ker(c^o_\Si)\rangle=\LR_\Si$.

In general, we notice that 
$$
E/G^o_\Xi=\BG^o_\Xi\xrar{\gamma_\Xi}E/G_\Xi=\BG_\Xi
\text{ and }
E/G^o_\Si=\BG^o_\Si\xrar{\gamma_\Si}E/G_\Si=\BG_\Si
$$
are principal $F$-bundles, and they induce isomorphism in 
the cohomology with rational coefficients. Then 
$$
H^*(\BG_\Xi;{\mbb Q})
\srel{\kern1ex \gamma_\Xi^*}{\cong} 
H^*(\BG^o_\Xi;{\mbb Q})
\cong 
\Sym^\bullet{\mbb Q}^{\Xi}
/\langle\Ker(u_{\Xi,\mbb Q})\rangle 
$$
and we are done again.
\end{proof}

Our main tool for computing the cohomology of $\crl X_\Si$ is 
the Leray-Hirsch theorem. For the confort of the reader, we 
recall this classical result:\smallskip

\nit\underbar{\sl Leray-Hirsch theorem} \cite[theorem 5.7.9]{sp} 
{\it Let $p:(\cal F,\dot{\cal F})\rar S$ be the total pair of 
a fibre-bundle pair with base $S$ and fibre pair $(F,\dot F)$, 
and let $R$ be a ring. Assume that $H^*(F,\dot F;R)$ is free 
over $R$ and that 
$\theta:H^*(F,\dot F;R)\rar H^*(\cal F,\dot{\cal F};R)$ is a 
cohomology extension of the fibre. Then the homomorphism
$$
\Theta^*:H^*(S;R)\otimes H^*(F,\dot F;R)\lar H^*(\cal F,\dot{\cal F};R),
$$ 
$$
u\otimes v\lmt p^*u\cup \theta(v)
$$
is an isomorphism of graded modules over $R$.}\smallskip

We should notice that in {\it loc.\ cit.} there is the additional 
assumption that $H^*(F,\dot F;R)$ is finitely generated over $R$. 
However, according to \cite[page 232]{hus}, this hypothesis can be 
omitted.

As an immediate consequence, we can determine the 
{\em ring structure} of such a fibre-bundle pair
$$
H^*(\cal F,\dot{\cal F};R)\cong 
H^*(S;R)\bigl[\theta(v)\mid v\in H^*(F,\dot F;R)\bigr],
$$
where in this writing the relations among the $\theta(v)$'s are 
understood.

We need also the following fact concerning the cohomology ring 
of complete and simplicial toric varieties:\smallskip

\nit\unl{\sl Result} \cite[theorem 10.8 and remark 10.9]{da}
{\it Let $X_\Si$ be a complete and simplicial toric variety. For 
each $\xi\in\Si(1)$ we denote $D_\xi$ the corresponding invariant 
divisor on $X_\Si$. Then the ring homomorphism 
$$
{\mbb Q}[\, f_\xi\mid \xi\in\Si(1)]\lar H^*(X_\Si;{\mbb Q})
$$ 
induced by the assignment
$$
f_\xi\lmt [D_\xi]\in H^2(X_\Si;{\mbb Q})
$$
is surjective. Its kernel is $\LR_\Si+\SR_\Si$, where 
$$
\begin{array}{l}
\LR_\Si \quad 
\text{\rm is the ideal of linear relations defined in 
lemma \ref{lm:bg}, and}
\\[1.5ex] 
\SR_\Si:=\bigl\langle
\prod_{\xi\in\pi}f_\xi\kern-1.5pt\mid\kern-1.5pt 
\pi\subset\Si(1)\text{ \rm is a primitive family}
\bigr\rangle
\;
\text{\rm is the Stanley-Reisner ideal.}
\end{array}
$$
If $X_\Si$ is smooth, that is $\Si$ is regular, the result holds 
with integral coefficients.}\smallskip

We shall use again a shorthand notation: we write $H^*(X_\Si)$ 
for the cohomology ring of $X_\Si$ with coefficients in $\mbb Z$ 
or $\mbb Q$, depending respectively on whether $\Si$ is regular 
or only simplicial. 

In view of the lemma \ref{lm:bg}, we deduce that
\begin{align}{\label{eq:coh-X}}
h_\Si:H^*_{G_\Si}\lar H^*(X_\Si),\quad \hat f_\xi\lmt [D_\xi]
\end{align}
is a surjective ring homomorphism and $\ker(h_\Si)=\SR_\Si$. 
Therefore the composition
$$
h:H^*_{G_\Xi}\xrar{\;b\;}H^*_{G_\Si}
\xrar{\;h_\Si\;} H^*(X_\Si)
$$
is still surjective with kernel 
$$
\Ker(h)=\langle f_j\mid j\in J\rangle+\SR_\Si
\subset H^*_{G_\Xi}.
$$

We can further say that the cohomology of $X_\Si$ is 
generated by the first Chern classes of the rank one 
sheaves $\cal O(D_\xi)\rar X_\Si$, $\xi\in\Si(1)$. 
Proposition \ref{prop:coord-f} gives a 
coordinate-free description for these sheaves as 
\begin{align}{\label{eq:Dxi}}
\cal O(D_\xi)
={\bigl(r_*\cal O_{U(\Si)\times W_0^\times}\bigr)}^{c_\Xi(f_\xi)},
\end{align}
for $r:U(\Si)\times W_0^\times\rar X_\Si$ the projection. 
Since $\Si$ is simplicial, $r$ is a geometric quotient, and we 
have a natural ring isomorphism 
$H^*(X_\Si)\srel{\cong}{\lar}H^*_{G_\Xi}\bigl(U(\Si)\bigr)$. It 
follows now from the second half of proposition \ref{prop:coord-f} 
that the diagram
\begin{align}{\label{cd4}}
\xymatrix{
H^*_{G_\Xi}(V)\ar[r]^-{\jmath^*}
&
H^*_{G_\Xi}\bigl(U(\Si)\times W_0^\times\bigr)
\\ 
H^*_{G_\Xi}\ar[u]^-{\cong}\ar[r]^-{h}
&
H^*(X_\Si)\ar[u]_-{\cong}
}
\end{align} 
commutes, so that $\jmath^*$ is surjective as $h$ is surjective. 
This fact, together with the long exact sequence in cohomology for 
the pair $\bigl(V,U(\Si)\bigr)$ imply that we have a short exact 
sequence
\begin{align}{\label{e-sqn1}}
0
\lar 
H^*_{G_\Xi}\bigl(V,U(\Si)\times W_0^\times\bigr)
\srel{\imath^*}{\lar} 
H^*_{G_\Xi}(V)
\srel{\jmath^*}{\lar}
H^*_{G_\Xi}\bigl(U(\Si)\times W_0^\times\bigr)
\lar 
0,
\end{align}
and $\imath^*$ and $\jmath^*$ are {\em ring homomorphisms}. 
We deduce now that 
$$
\imath^*H^*_{G_\Xi}\bigl(V,U(\Si)\times W_0^\times\bigr)
\subset H^*_{G_\Xi}(V)
$$
is an ideal. In view of the previous discission, this ideal 
is precisely
\begin{align}{\label{eq:ideal}}
\begin{array}{ll}
\imath^*H^*_{G_\Xi}\bigl(V,U(\Si)\times W_0^\times\bigr)
&
=
\bigl\langle
\prod_{\xi\in\pi}\hat f_\xi \mid \pi\text{ primitive family}
\bigr\rangle
+\langle \hat f_j \mid j\in J \rangle
\\[1.5ex]
&
\disp
=
\biggl\langle e_{G_\Xi}(W_\pi)\;\biggl|\;
\pi\text{ primitive family and }
W_\pi:=\bigoplus_{\beta\subset\pi}W_\beta 
\biggr.\biggr\rangle
\\[1.5ex]
&
\disp\kern3ex
+\bigl\langle e_{G_\Xi}(V_j)\ \bigl|\ j\in J\bigr\rangle.
\end{array}
\end{align}
In this formula, $e_{G_\Xi}(\ \cdot\ )$ stays for the equivariant 
Euler class of the correponding vector space, and we have used that 
$$
e_{G_\Xi}(W_\pi)
=
\prod_{\beta\subset\pi}
\biggl[
\prod_{\alpha\in b^{-1}(\beta)}\kern-1ex e_{G_\Xi}(V_\alpha)
\biggr]
=
\prod_{\beta\subset\pi}
\biggl[
\prod_{\alpha\in b^{-1}(\beta)}
\biggl(\prod_{c_\Xi(f_\xi)=\alpha}\kern-1ex\hat f_\xi\biggr)
\biggr]
\srel{\ref{prop:prim}}{=}
\prod_{\xi\in \pi}\hat f_\xi.
$$
The conclusion of our discussion can now be formulated as

\begin{corollary}{\label{cor:gener}}
The classes ${(\imath^*)}^{-1}e_{G_\Xi}(W_\pi)$, where 
$\pi\subset\Si(1)$ ranges over the primitive families, 
and ${(\imath^*)}^{-1}e_{G_\Xi}(V_j)$, $j\in J$, generate 
$H^*_{G_\Xi}\bigl(V,U(\Si)\bigr)$ as a $H^*_{G_\Xi}$-algebra.
\end{corollary}

We should point out that is not so easy to compute the 
explicit formulae for these generators.

\begin{proposition}{\label{prop:lh}}
The Leray-Hirsch theorem does apply to the fibration
$$
\begin{array}{r}
X_\Si\hra\tcx_\Si
\\ 
\dar\kern1.4ex 
\\ 
S.\kern.8ex
\end{array}
$$
\end{proposition}

\begin{proof} We must show that there does exist a cohomology 
extension of the fibre $X_\Si$. We have just seen that the 
cohomology of $X_\Si$ is generated by the first Chern classes 
of the sheaves $\cal O(D_\xi)$, $\xi\in \Si(1)$, and that they 
can be described by \eqref{eq:Dxi}. In the relative setting, 
for each $\xi\in\Si(1)$, we have defined in proposition 
\ref{prop:lxi} the rank one sheaf 
\begin{align*}
\tcl_\xi:={\bigl({r}_*\cal O_{\tcu(\Si)}\bigr)}^{c_\Xi(f_\xi)},
\end{align*}
where $r:\crl U(\Si)\rar \tcx_\Si$ is the projection, 
whose restriction to the fibres of $\tcx_\Si\rar S$ is 
isomorphic to $\cal O(D_\xi)$. 

We claim that the sheaves $\tcl_\xi$ possess a first Chern 
class. Indeed, we have proved in proposition \ref{prop:lxi} 
that $\tcl_\xi$ is locally free outside of a subvariety 
$\tcx_{\Si, \rm bad}\hra\tcx_\Si$ of codimension at least 
two. Then the first Chern class of 
$\tcl_\xi|_{\tcx_\Si\sm \tcx_{\Si,\rm bad}}$ is a well-defined 
element of $H^2(\tcx_\Si\sm \tcx_{\Si,\rm bad};{\mbb Z})$, 
and this group is isomorphic to $H^2(\tcx_\Si;{\mbb Z})$ since 
the bad set has {\em complex} codimension at least two in 
$\tcx_\Si$.
\end{proof}

\begin{corollary} 
The cohomology ring of $\tcx_\Si$ is
$$
H^*(\tcx_\Si)\cong H^*(S)
\bigl[\, l_\xi \mid \xi\in\Si(1)\bigr],
\quad\text{with}\quad 
l_\xi:=c_1(\tcl_\xi)\in H^2(\crl X_\Si;{\mbb Z}).
$$
In particular, the ring homomorphism
$$
h_S: H^*_{G_\Xi}(S)=H^*(S)\otimes H^*_{G_\Xi}\lar H^*(\tcx_\Si),
$$ 
$$
u\otimes \hat f_\rho\lmt \tld\phi^*u\cup l_\rho,\;\forall\ \rho\in\Xi,
$$
is surjective, where $\tld\phi:\tcx_\Si\rar S$ is the projection.
\end{corollary}

Before entering into the proof, we should clarify the 
notation: for $\rho=j\in J\subset\Xi$, we know from 
proposition \ref{prop:lxi} that 
$\tcl_j\cong\tld\phi^*\euf V_j^\vee$, 
and therefore $l_j=-\tld\phi^*[\euf V_j]$.

\begin{proof}
All we need to prove is that the homomorphism $h_S$ 
is well-defined. Since $H^*_{G_\Xi}$  is the quotient 
$H^*_\cxi\bigl/\LR_\Xi\bigr.$, we must check that as soon 
as the classes ${\{\hat f_\rho\}}_\rho$ satisfy a linear 
relation, the classes ${\{l_\rho\}}_\rho$ fulfill the same 
one. 

So, let us consider integers ${\{a_\rho\}}_\rho$ having the 
property that $\sum_\rho a_\rho\hat f_\rho=0$. This is 
equivalent saying that $\sum_\rho a_\rho c_\Xi(f_\rho)
=\sum_\rho c_\Xi(a_\rho f_\rho)=0\in \cal X^\star(G_\Xi)$. 
We are going to prove that 
$\bigotimes_\rho \tcl_\rho^{\otimes a_\rho}
\cong \cal O_{\tcx_\Si}$ over 
$\tcx_\Si\sm \tcx_{\Si,\rm bad}$, and the conclusion follows 
then by a codimension two argument. 

Let us consider an open subset 
$O\subset \tcx_\Si\sm \tcx_{\Si,\rm bad}$, and 
non-vanishing sections $s_\rho\in\Gamma(O,\tcl_\rho)
={\Gamma\bigl(r^{-1}O,\cal O_{\tcu(\Si)}\bigr)}^{c_\Xi(f_\rho)}$. 
Then 
$$
s:=\bigotimes_\rho s_\rho^{\otimes a_\rho}
\in
{\Gamma\bigl( r^{-1}O,\cal O_{\tcu(\Si)}\bigr)}
^{\sum_\rho c_\Xi(a_\rho f_\rho)=0}
=
{\Gamma\bigl( r^{-1}O,\cal O_{\tcu(\Si)}\bigr)}^{G_\Xi}
=
\Gamma\bigl( O, \cal O_{\tcx_\Si}\bigr).
$$
The last equality is implied by the fact that $\tcx_\Si$ 
is a geometric quotient of $\tcu(\Si)$ for the $G_\Xi$-action.
\end{proof}

In order to make computations, one still needs to determine 
the relations between the classes $l_\rho$, $\rho\in\Xi$. 
Proposition \ref{prop:lxi} implies that the diagram 
\begin{align}{\label{cd5}}
\xymatrix{
H^*_{G_\Xi}(\euf V)\ar[r]^-{\jmath^*_S}
&
H^*_{G_\Xi}\bigl(\tcu(\Si)\bigr)
\\ 
H^*_{G_\Xi}(S)\ar[u]^-\cong\ar[r]^-{h_S}
&
H^*(\tcx_\Si)\ar[u]_-\cong
}
\end{align}
commutes. Since $h_S$ is surjective, $\jmath^*_S$ is surjective 
too, and it follows from the long exact sequence in equivariant 
cohomology for the pair $\bigl(\euf V,\tcu(\Si)\bigr)$ 
that we have the short exact sequence
\begin{align}{\label{e-sqn2}}
0
\lar 
H^*_{G_\Xi}\bigl(\euf V,\tcu(\Si)\bigr)
\srel{\imath^*_S}{\lar}
H^*_{G_\Xi}(\euf V)
\srel{\jmath^*_S}{\lar}
H^*_{G_\Xi}\bigl(\tcu(\Si)\bigr)
\lar 
0,
\end{align}
and again $\imath^*_S$ and $\jmath^*_S$ are ring homomorphisms. 
Therefore 
$$
\imath^*_S H^*_{G_\Xi}\bigl(\euf V,\tcu(\Si)\bigr)\subset 
H^*_{G_\Xi}(\euf V)
$$
is an ideal, and our next task is to compute it. 

\begin{proposition}{\label{prop:rel-lh}}
The Leray-Hirsch theorem does apply for the fibre bundle pair 
$$
\begin{array}{r}
\bigl(E\times_{G_\Xi} V,E\times_{G_\Xi}
 (U(\Si)\times W_0^\times)\bigr)
\hra
\bigl(E\times_{G_\Xi}\euf V,E\times_{G_\Xi}\tcu(\Si)\bigr)
\\ 
\dar\kern13.5ex 
\\ 
S.\kern12.8ex 
\end{array}
$$
\end{proposition}

\begin{proof}
We have already described in the corollary \ref{cor:gener} the 
generators for the cohomology of the fibre, and we must prove 
that they extend to classes in 
$H^*_{G_\Xi}\bigl(\euf V,\tcu(\Si)\bigr)$. 
First of all, we should find the candidates for the cohomology 
extension. For any $\Si$-primitive family $\pi$, we define 
$$
\euf W_\pi:=
\bigoplus_{\beta\subset\pi}\unl{\mbb C}_\beta\otimes\euf W_\beta
$$
and we consider the projection
$$
\pr_\pi:\euf V\lar \euf W_\pi.
$$
Denoting by $\pr_S:\euf V\rar S$ the projection, we define the 
section
$$
s_\pi:\euf V\lar \pr_S^*\euf W_\pi,\quad 
v\lmt \bigl(v,\pr_\pi(v)\bigr)
$$
which is obviously $G_\Xi$-equivariant, and it vanishes nowhere 
on $\tcu(\Si)\subset\euf V$. Consequently, the restriction of 
the equivariant Euler class of $\pr_S^*\euf W_\pi$ to 
$\tcu(\Si)$ vanishes
$$
\jmath^*_S\bigl( e_{G_\Xi}(\pr_S^*\euf W_\pi)\bigr)=0
\in H^*_{G_\Xi}\bigl(\tcu(\Si)\bigr).
$$
The exactness of the sequence \eqref{e-sqn2} implies now that 
$$
e_{G_\Xi}(\pr_S^*\euf W_\pi)\in 
\imath^*_S H^*_{G_\Xi}\bigl(\euf V,\tcu(\Si)\bigr),\quad 
\text{ for any primitive family }\pi.
$$

For every $j\in J$, we have the $G_\Xi$-equivariant section 
$$
s_j:\euf V\rar \pr_S^*\bigl(\unl{\mbb C}_j\otimes\euf V_j\bigr),
\quad  v\lmt (v,v_j).
$$
Again, this section vanishes nowhere on $\tcu(\Si)$, and 
therefore the restriction of the equivariant Euler class 
$e_{G_\Xi}\bigl(\pr_S^*(\unl{\mbb C}_j\otimes\euf V_j)\bigr)$ 
to $\tcu(\Si)$ is zero. The exactness of \eqref{e-sqn2} 
implies again that 
$$
e_{G_\Xi}\bigl(\pr_S^*(\unl{\mbb C}_j\otimes\euf V_j)\bigr)
\in 
\imath^*_S H^*_{G_\Xi}\bigl(\euf V,\tcu(\Si)\bigr),\quad 
\text{ for any }j\in J.
$$

It is now clear, by corollary \ref{cor:gener}, that the 
cohomology classes 
${(\imath^*_S)}^{-1}\bigl( e_{G_\Xi}(\pr_S^*\euf W_\pi)\bigr)$ 
and 
$e_{G_\Xi}\bigl(\pr_S^*(\unl{\mbb C}_j\otimes\euf V_j)\bigr)$ 
define, in fact generate, a cohomology extension of the 
fibre bundle pair.
\end{proof}

\begin{corollary}{\label{gener}}
The classes 
${(\imath^*_S)}^{-1}\bigl(e_{G_\Xi}(\pr_S^*\euf W_\pi)\bigr)$, 
where $\pi$ ranges over the $\Si$-primi\-tive families, and 
${(\imath^*_S)}^{-1}
\bigl(e_{G_\Xi}\bigl(\pr_S^*(\unl{\mbb C}_j\otimes
\euf V_j)\bigr)\bigr)$, 
$j\in J$, generate 
$H^*_{G_\Xi}\bigl(\euf V,\tcu(\Si)\bigr)$ as a 
$H^*_{G_\Xi}(S)$-alge\-bra.
\end{corollary}

Since the cohomology rings $H^*_{G_\Xi}(\euf V)$ and 
$H^*_{G_\Xi}(S)$ are naturally isomorphic, we shall identify 
$e_{G_\Xi}(\pr_S^*\euf W_\pi)$ and 
$e_{G_\Xi}\bigl(\pr_S^*(\unl{\mbb C}_j\otimes\euf V_j)\bigr)$ 
respectively with $e_{G_\Xi}(\euf V_\pi)$ and 
$e_{G_\Xi}\bigl(\unl{\mbb C}_j\otimes\euf V_j\bigr)$. The 
explicit formulae for them are  
$$
e_{G_\Xi}(\euf W_\pi)
=\prod_{b(\alpha)\subset\pi}\kern-1ex
e_{G_\Xi}\bigl(\unl{\mbb C}_\alpha\otimes\euf V_\alpha\bigr)
=
\prod_{b(\alpha)\subset\pi}\kern-1ex
\biggl(
\alpha^{r_\alpha}
+c_1(\euf V_\alpha)\alpha^{r_\alpha-1}
+{\dots}
+c_{r_\alpha}(\euf V_\alpha)
\biggr),
$$ 
$$
e_{G_\Xi}\bigl(\unl{\mbb C}_j\otimes\euf V_j\bigr)
=c_1(\euf V_j)+c_\Xi(f_j)=c_1(\euf V_j)+\hat f_j.
$$

\begin{definition}{\label{def:sr}}
We define the ideal of linear relations and the 
Stanley-Reisner ideal, both relative to $S$, 
respectively by
$$
\begin{array}{l}
\cLR_\Xi:=H^*(S)\otimes \LR_\Xi,
\quad
\cLR_\Si:=H^*(S)\otimes \LR_\Si,
\\[1.5ex]
\cSR_\Xi:=
\bigl\langle 
e_{G_\Xi}(\euf W_\pi)\mid\pi\text{ is a $\Si$-primitive family}
\bigr\rangle
\subset H^*_{G_\Xi}(S)=H^*_{G_\Xi}\otimes H^*(S),
\\[1.5ex]
\cSR_\Si:=
\bigl\langle 
e_{G_\Si}(\euf W_\pi)\mid\pi\text{ is a $\Si$-primitive family}
\bigr\rangle
\subset H^*_{G_\Si}(S)=H^*_{G_\Si}\otimes H^*(S).
\end{array}
$$
By abuse of notation, in the next theorem, we shall still 
write $\cSR_\Xi$ and $\cSR_\Si$ for the inverse images of 
these ideals under the epimorphisms 
$H^*_{\cxi}\otimes H^*(S)\rar H^*_{G_\Xi}\otimes H^*(S)$ 
and 
$H^*_{\csi}\otimes H^*(S)\rar H^*_{G_\Si}\otimes H^*(S)$
respectively.
\end{definition}

Now we are finally in position to compute the cohomology 
rings of $\tcx_\Si$ and $\crl X_\Si$.

\begin{theorem}{\label{the-cohom}}
Let $\Si$ be a complete and simplicial fan, and consider 
$\euf V\rar S$ a vector bundle. 

\nit {\rm (i)} The cohomology ring of the family 
$\tcx_\Si\rar S$ is 
$$
\begin{array}{ll}
H^*(\tcx_\Si;\mbb Q)
&
\cong 
\raise.25ex\hbox{$H^*_{G_\Xi}(S;\mbb Q)$}
\bigl/
\lower.25ex
\hbox{$\cSR_\Xi+\langle c_1(\euf V_j)+\hat f_j, j\in J\rangle$}
\bigr.
\\[1.5ex]
&
\cong 
\raise.25ex\hbox{$H^*(S;\mbb Q)[f_\rho\mid\rho\in\Xi]$}
\bigl/
\lower.25ex\hbox{$\cLR_\Xi+\cSR_\Xi+
\langle c_1(\euf V_j)+f_j, j\in J\rangle$}
\bigr.
\end{array}
$$

\nit {\rm (ii)} The cohomology ring of the family 
$\crl X_\Si\rar S$ is 
$$
\begin{array}{ll}
H^*(\tcx_\Si;\mbb Q)
&
\cong 
\raise.25ex\hbox{$H^*_{G_\Si}(S;\mbb Q)$}
\bigl/
\lower.25ex
\hbox{$\cSR_\Si$}
\bigr.
\cong 
\raise.25ex\hbox{$H^*(S;\mbb Q)[f_\xi\mid\xi\in\Si]$}
\bigl/
\lower.25ex\hbox{$\cLR_\Si+\cSR_\Si$}.
\bigr.
\end{array}
$$

\nit If $\Si$ is regular, both isomorphisms hold with 
$\mbb Z$-coefficients.
\end{theorem} 

\begin{proof}
(i) The key for the proof is the exact sequence \eqref{e-sqn2}: 
it implies that, as a ring, 
$$
H^*(\tcx_\Si)
\cong 
\raise.25ex\hbox{$H^*_{G_\Xi}(\euf V)$}\bigl/
\lower.25ex\hbox{$\imath^*_S H^*_{G_\Xi}
\bigl(\euf V,\tcu(\Si)\bigr)$}.
\bigr.
$$
We should point out that since $\jmath_S^*$ is a ring 
homomorphism, 
$\imath^*_S H^*_{G_\Xi}\bigl(\euf V,\tcu(\Si)\bigr)\subset 
H^*_{G_\Xi}(\euf V)$ is an ideal, and the r\^ole of the 
previous discussion was to determine this ideal. Corollary 
\ref{gener} implies that 
$$
\begin{array}{ll}
\imath_S^* H^*_{G_\Xi}\bigl(\euf V,\tcu(\Si)\bigr)
&\disp
\cong
\sum_{\pi\ \rm primitive}\kern-1.5ex
H^*_{G_\Xi}(S)\cdot e_{G_\Xi}(\euf W_\pi)
+\langle c_1(\euf V_j)+\hat f_j, j\in J\rangle
\\[1.5ex]
&
=\cSR_\Xi+ \langle c_1(\euf V_j)+\hat f_j, j\in J\rangle.
\end{array}
$$
The first isomorphism in the statement of theorem becomes 
now clear. The second isomorphism is an immediate consequence 
of the fact that $H^*_{G_\Xi}(S)=H^*(S)\otimes H^*_{G_\Xi}$, 
and of the lemma \ref{lm:bg} which states that 
$H^*_{G_\Xi}\cong H^*_{\cxi}\bigl/\LR_\Xi\bigr.$.

\nit (ii) The structure of the cohomology ring of the 
family $\crl X_\Si\rar S$ is obtained by replacing 
$\euf V$ with $\euf W_\Si$, or equivalently in the 
particular case when $\Si(1)=\Xi$. 
\end{proof}

What we should remark about this result is that the cohomology 
rings of the families $\tcx_\Si$ and $\crl X_\Si$ are not isomorphic 
in general, so that these families are really distinct despite the 
fact that both are $X_\Si$-fibrations over $S$. The reason for this 
phenomenon is that the isomorphism found in proposition 
\ref{prop:inv-quot} was not natural.


\section{A special case}

Our primary motivation for studying families of toric 
varieties was their occurence in the context of the gauged 
sigma models. In this section we are going to put the 
phenomenon which occurs there in our present framework. 
In \cite{ha} and \cite{ot}, the starting point was to consider 
a regular, projective fan $\Si\subset N_{\mbb R}$ and its 
associated toric variety $X_\Sigma$. The first question one 
has to answer in order to compute enumerative invariants of 
$X_\Si$ is to describe the space of morphisms, with 
fixed degree, from a smooth and projective curve $C$ into 
$X_\Sigma$. It turns out that when this degree is sufficiently 
large, the compactification of this space of morphisms has the 
structure of a fibre bundle over the product of $b_2(X_\Si)$ 
copies of the Jacobian of $C$. In our language, it is a family 
of toric varieties parameterized by the base 
$S={\rm Jac}(C)^{b_2(X_\Si)}$. 

The fibres of this family are toric varieties which `look very 
much the same' as the starting variety $X_\Si$. It is the 
meaning of this sentence which we intend to explain here. So let 
us consider a finite subset $\Xi\subset N$ having the properties 
enumerated in section \ref{sct:rmk}. The diagram \eqref{cd1} 
provides us with the exact sequence
$$
0\lar M\srel{a}{\lar}{\mbb Z}^\Xi\srel{c}{\lar}A:=A_\Xi\lar 0, 
$$
and therefore we can decompose 
$$
{\mbb Z}^\Xi=\bigoplus_{\xi\in\Xi}{\mbb Z}f_\xi
=\bigoplus_{\alpha\in A}
\biggl(
\bigoplus_{f_\xi\in c^{-1}(\alpha)}{\mbb Z}f_\xi
\biggr)
=:\bigoplus_{\alpha\in A} L_\alpha.
$$
Each of the $L_\alpha$'s are lattices, since they are torsion 
free $\mbb Z$-modules. Let us choose arbitrary, strictly positive 
integers $\{n_\alpha\}_{\alpha\in A}$ and define 
$$
L':=\bigoplus_{\alpha\in A} L_\alpha^{\oplus n_\alpha}.
$$
This amounts saying that for every $\alpha\in A$ we reproduce  
a number of $n_\alpha$ times each vector $f_\xi\in c^{-1}(\alpha)$. 
For fixed $\alpha\in A$ and $f_\xi\in c^{-1}(\alpha)$, we denote 
by $\{f_{\xi 1},\dots,f_{\xi n_\alpha} \}$ these vectors, so we 
can write 
$$
L'=\bigoplus_{\alpha\in A}
\biggl(
\bigoplus_{f_\xi\in c^{-1}(\alpha)}
\bigl( 
{\mbb Z}f_{\xi 1}\oplus{\dots}\oplus{\mbb Z}f_{\xi n_\alpha}
\bigr)
\biggr).
$$
It is now clear that there is a natural surjective homomorphism 
$$
c':L'\lar A,\quad f_{\xi \bfd}\lmt c(f_\xi)\in A.
$$
Let us denote by $M':=\ker(c')$ and by $a':M'\rar L'$ the inclusion. 
We notice that $M'$ is a lattice, since without torsion, whose rank 
is
$$
{\rm rk}(M')
={\rm rk}(M)+\sum_{\alpha\in A}(n_\alpha-1){\rm rk}(L_\alpha).
$$
The lattice $L'$ `looks very much like' ${\mbb Z}^\Xi$ in the 
following sense: for every $\alpha\in A$ and 
$f_\xi\in c^{-1}(\alpha)$, and for any choice of 
$\bfd_\xi\in\{1,\dots,n_\alpha\}$, we find the 
monomorphism 
$$
j_\bfd:{\mbb Z}^\Xi\lar L',\quad f_\xi\lmt f_{\xi \bfd_\xi}.
$$
There is also a canonical epimorphism 
$$
+:L'=\bigoplus_{\alpha\in A} L_\alpha^{n_\alpha} 
\lar 
\bigoplus_{\alpha\in A} L_\alpha={\mbb Z}^\Xi,
$$
induced by the natural additions 
$+_\alpha:L_\alpha^{n_\alpha}\rar L_\alpha$, $\alpha\in A$. 
Let us notice now that for all choices of $j_\bfd$ the diagram
\begin{align*}
\xymatrix{
0\ar[r]
& 
M\ar[d]\ar[r]^-a
&
{\mbb Z}^\Xi \ar@/^/[d]^-{j_\bfd}\ar[r]^-c
&
A\ar[r]\ar@{=}[d]
&
0
\\ 
0\ar[r]
& 
M'\ar[r]^-{a'}
&
L' \ar@/^/[u]^-{+}\ar[r]^-{c'}
&
A\ar[r]
&
0
}
\end{align*}
commutes. By dualizing, we get for $N':=\Hom_{\mbb Z}(M',{\mbb Z})$ 
the commutative diagram
$$
\xymatrix{
{(L')}^\vee
=\bigoplus_{\xi,\bfd_\xi}{\mbb Z}f_{\xi \bfd_\xi}^\vee
\ar[r]^-{(a')^\vee}\ar@{->>}[d]^-{j_\bfd^\vee}
&
N'\ar[d]^-{\nu_\bfd} & 
\\ 
{\mbb Z}_\xi=\bigoplus_{\xi\in\Xi}{\mbb Z}f_\xi^\vee
\ar[r]^-{a^\vee}
&
N\ar[r]
&
0.
}
$$
and observe that the homomorphism $\nu_\bfd$ is surjective. 
We define now the set 
$\Xi':=\bigl\{(a')^\vee(f_{\xi \bfd}^\vee)\bigr\}\subset N'$, 
where $f_{\xi \bfd}^\vee$ stays for the canonical generator of 
$\Hom_{\mbb Z}({\mbb Z}f_{\xi \bfd},{\mbb Z})$. 
Our goal is to prove the following\medskip 

\begin{lemma}
The set $\Xi'\subset N'$ fulfills the conditions (i)-(iii) 
and $L'\cong {\mbb Z}^{\Xi'}$.
\end{lemma}

\begin{proof} We claim that $0\not\in \Xi'$: otherwise we find 
an $f_{\xi_0 \bfd_{\xi_0}}$ having the property that all 
vectors in $M'$ have the corresponding entry zero in $L'$. 
We choose now an inclusion $j_\bfd$ such that 
$f_{\xi_0}\mt f_{\xi_0 \bfd_{\xi_0}}$. Then it follows from 
the diagram above that 
$$
\xi_0=a^\vee(f_{\xi_0}^\vee)
=a^\vee\bigl(j_\bfd^\vee(f_{\xi_0 \bfd_{\xi_0}}^\vee)\bigr)=0, 
$$
and this is a contradiction. 

We are going to check now the second property, that 
each $\xi'=(a')^\vee(f_{\xi \bfd_\xi})$ is the generator 
over $\mbb Z_{\geq 0}$ of the semi-group 
${\mbb R}_{\geq 0}\xi'\cap N'$. We chose an appropriate 
inclusion $j_\bfd$ as above. Then we know that 
$\xi=\nu_\bfd(\xi')=a^\vee(f_\xi^\vee)$ is the generator 
of the semi-group ${\mbb Z}_{\geq 0}\xi\cap N$ and 
therefore $\xi'$ has the same property in $N'$. 

It remains to check that 
$\sum_{\xi'\in \Xi'}{\mbb R}_{\geq 0}\xi'=N'_{\mbb R}$. The left 
hand side is a convex cone in $N'_{\mbb R}$, and the condition is 
equivalent asking that there is no hyperplane in $N'_{\mbb R}$ 
which leaves this cone on one of its sides. In other words, it 
is enough checking that there is no $m'\in M'$ such that 
$\langle m', \xi'\rangle\geq 0$, for all $\xi'\in \Xi'$. Let us 
assume the contrary that such an $m'\in M'\sm\{0\}$ does exist. 
Then $(c\circ+\circ a')(m')=0\in A$, and consequently 
$(+\circ a')(m')=a(m)$ for some $m\in M$. This $m$ is certainly 
not zero, because $a(m')$ has posivive entries and has at least 
one strictly positive entry. Then we have found an element 
$m\in M\sm \{0\}$ which leaves all the $\xi\in\Xi$ on one of its 
sides. This is a contadiction with the fact that $\Xi\subset N$ 
does obey (iii). 

We must still check a tacit assumption, namely that the $\xi'$'s 
are indeed pairwise distinct. Let us assume the contrary, 
that there exist $f_{\xi_1\bfd_1}$ and $f_{\xi_2\bfd_2}$ which 
determine the same element in $N'$. This means that for any 
$m'\in M'$, the element $a'(m')\in L'$ has the 
${\xi_1\raise.3ex\hbox{\bfd}_1}$ 
and the ${\xi_2\raise.3ex\hbox{\bfd}_2}$ coordinates equal. 
But if follows from the definition of $M'$ that 
$$
M'=\ker(c')=\ker(c\circ +)
=\left\{
\sum a_{\xi\bfd_\xi}f_{\xi\bfd_\xi}
\left| 
\begin{array}{l}
\exists m\in M\text{ s.t. }\forall\alpha\in A,\\  
\disp
\sum_{
\xi\in c^{-1}(\alpha),\; 
1\leq\bfd_\xi\leq n_\alpha
}
\kern-3ex a_{\xi\bfd_\xi}=\langle m,\xi\rangle
\end{array}
\right.
\right\}.
$$
We distinguish between two possibilities: the first one is when 
$\xi_1=\xi_2=\xi$. We can choose in this case the vector whose 
$\xi\raise.3ex\hbox{\bfd}_1$ and $\xi\raise.3ex\hbox{\bfd}_2$ 
entries are $1$ and $-1$ respectively, and zero in rest. It 
clearly belongs to $\ker(c')=M'$, but however the entries of 
interest are different. 
The second case is when $\xi_1\neq \xi_2$: then we can choose 
a monomorphism $j_\bfd:{\mbb Z}^{\Si(1)}\rar L'$ which sends 
$f_{\xi_1}\mt f_{\xi_1\bfd_1}$ and $f_{\xi_2}\mt f_{\xi_2\bfd_2}$. 
We deduce then that 
$\xi_1=a^\vee(f_{\xi_1})=a^\vee(f_{\xi_2})=\xi_2$, which is again 
a contradiction. Therefore the vectors $\xi'\in N'$ defined above 
are pairwise distinct.

The last property in the lemma follows now from the definions of 
$N'$ and $\Xi'$.
\end{proof}

We can finally explain the situation which appears in 
\cite{ha, ot}: $\Si\in N_{\mbb R}$ is a projective and regular 
fan, and $\Xi=\Si(1)$. One defines a set $\Xi'\subset N'$ as 
above, and the space of morphisms with fixed degree from a 
curve $C$ into $X_\Si$ turns out to be a toric fibration over 
${\rm Jac}(C)^{b_2(X_\Si)}$. The fibres are toric varieties 
associated to the fan $\Si'\subset N'_{\mbb R}$ whose 
primitive families are 
$$
\bigl\{
(a')^\vee\bigl(f_{\xi\bfd_\xi}^\vee\bigr)
\mid \xi\in \pi,\ \bfd_\xi\in\{1,\dots,n_{c(f_\xi)}\}
\bigr\},
$$ 
with $\pi\subset\Xi$ a $\Si$-primitive family. The maximal 
cones of this fan can be described as follows: for each cone 
$\si\in\Si_{\rm max}$, choose $k_\xi\in\{1,\dots,n_{c(f_\xi)}\}$, 
for all $\xi\not\in\si(1)$; define then 
$\si'\in\Si'_{\rm max}$ to be the cone whose $1$-skeleton is 
\begin{align}{\label{sip}}
\si'(1)=
\left\{
(a')^\vee\bigl(f_{\xi\bfd_\xi}\bigr)\ 
\left|\ 
\begin{array}{l}
\bfd_\xi\in\{1,\dots,n_{c(f_\xi)}\},
\quad\forall\ \xi\in\si(1)
\\[1ex] 
\bfd_\xi\in\{1,\dots,n_{c(f_\xi)}\}\sm\{k_\xi\},
\quad\forall\ \xi\not\in\si(1)
\end{array}
\right.
\right\}.
\end{align}

\begin{lemma}
The fan $\Si'$ is projective and regular. 
\end{lemma}

\begin{proof}
Since $\Si$ is regular, $X_\Si$ is isomorphic to the geometric 
quotient $U(\Si)/G$, and $U(\Si)\rar X_\Si$ is a principal 
$G$-bundle, where $G=\Hom_{\mbb Z}(A,{\mbb C}^*)$. Moreover, 
since $\Si$ is projective, there exists an ample class 
$\chi\in A=A^1(X_\Si)$, and $X_\Si$ can also be described as 
the geometric invariant quotient ${\mbb C}^\Xi\invq G$ for 
the linearization 
$$
G\times\bigl( {\mbb C}^\Xi\times {\mbb C}\bigr)\lar 
{\mbb C}^\Xi\times {\mbb C},\quad 
g\times(z,\veps)=(g\times z,\chi(g)\veps).
$$

Since $\sum_{\xi'\in\Xi'}{\mbb R}_{\geq 0}\xi'=N'$, we 
deduce that $\Si'$ is complete. Cox's result says that 
$X_{\Si'}$ is (canonically) isomorphic to the categorical 
quotient $U'(\Si')/G$, and 
\begin{align}{\label{usip}}
U'(\Si')={\mbb C}^{\Xi'}\sm Z_{\Si'}
=\bigcup_{\si'\in\Si'_{\rm max}} U'_{\si'}\quad{\rm with}\quad 
U'_{\si'}=
\biggl\{
\prod_{\xi\not\in\si(1)} Z_{\xi k_\xi}\neq 0.
\biggr\}.
\end{align}
We linearize the $G$-action on ${\mbb C}^{\Xi'}$ in the trivial 
line bundle, using the character $\chi$, and we want to prove 
that 
\begin{align}{\label{incl}}
U'(\Si')\subseteq \bigl({\mbb C}^{\Xi'}\bigr)^{{\rm ss}(\chi)}.
\end{align}
We will deduce then that $X_{\Si'}$ is quasi-projective and 
complete, and is consequently projective. {\it A posteriori} 
we find that the inclusion above is in fact an equality.

Let us consider a maximal cone $\si'$ of $\Si'$, and choose an 
inclusion $j_\bfd:{\mbb C}^\Xi\rar {\mbb C}^{\Xi'}$ having the 
property that $Z_\xi\mt Z_{\xi k_\xi}$, for all $\xi\not\in\si(1)$. 
We denote by  $p_\bfd:{\mbb C}^{\Xi'}\rar {\mbb C}^\Xi$ 
the corresponding projection, and we notice that both $j_\bfd$ 
and $p_\bfd$ are $G$-equivariant. Moreover, equality \eqref{usip} 
implies that $U_\si=p_\bfd(U'_{\si'})$. 

If $z'\in U'(\Si')$, there is a maximal cone $\si'\in\Si'$ such 
that $z'\in U'_{\si'}$. Then $z=p_\bfd(z')\in U_\si\subset 
U(\Si)$, and we find an integer $\kappa\geq 1$ and $s\in 
{{\mbb C}[{\mbb C}^\Xi]}^{\chi^\kappa}$ such that $s(z)\neq 0$. 
Then $(p_\bfd)^*s$ is a $\chi^\kappa$-invariant section on 
${\mbb C}^{\Xi'}$ which does not vanish at $z'$. The inclusion 
\eqref{incl} is proved now. 

We must still show that $\Si$ is regular. The first thing to 
notice is that $U'(\Si')\rar X_{\Si'}$ is actually a geometric 
quotient. It is enough to check that all $G$-orbits are closed 
in $U'_{\si'}$, for every $\si'\in\Si'_{\rm max}$. This follows 
from the following facts: first, since the quotient 
$U(\Si)\rar X_\Si$ is geometric, all the $G$-orbits are closed 
in $U_\si$, for all $\si\in\Si_{\rm max}$; secondly, the image 
of the inclusion $j_\bfd:{\mbb C}^\Xi\rar {\mbb C}^{\Xi'}$ 
defined above has closed image and 
$U_\si=j_\bfd^{-1}(U'_{\si'})$. Then its restiction to $U_\si$ 
still has closed image in $U'_{\si'}$ and is $G$-equivariant.  
Therefore the $G$-orbits in $U'_{\si'}$ are closed, hence the 
quotient $U'(\Si')\rar X_{\Si'}$ is geometric, as stated. 
Using again Cox's theorem \cite[theorem 2.1]{cox}, we deduce 
that $\Si'$ is simplicial.

Moreover, $\Si$ being regular, $G$ acts freely on each affine 
piece $U_\si$. Using the inclusion $j_\bfd$, we deduce that $G$ 
acts freely on $U'_{\si'}$ too, so that $X_{\Si'}$ is smooth. 
We conclude now from \cite[page 29]{fu} that $\Si'$ is regular.
\end{proof}

The statement, at the beginning of this section, that the 
toric variety $X_{\Si'}$ `looks the same' as $X_\Si$ means, 
after this discussion, that through each point of $X_{\Si'}$ 
passes a 
$\sum_{\alpha\in A}(n_\alpha-1){\rm rk}(L_\alpha)$-dimensional 
algebraic family of embedded copies of $X_\Si$. 

The simplest exemple is when $X_\Si={\mbb P}^{r-1}$ and 
$X_{\Si'}={\mbb P}^{nr-1}$. In this case through each point 
of ${\mbb P}^{nr-1}$ passes a family of linearly embedded 
copies of ${\mbb P}^{r-1}$'s, parameterized by 
${\mbb P}^{(n-1)r}$.


\section*{Appendix}

For proving proposition 4.6, D.~Cox makes use in \cite[page 43]{cox} 
of the following\smallskip 

\nit{\sl Lemma}\quad 
{\it Let $\Si\subset N_{\mbb R}$ be a complete fan 
and consider $m\in M=N^\vee$ such that the corresponding hyperplane 
${\cal H}_m:=\{\langle m,\ \cdot\ \rangle=0\}\subset N_{\mbb R}$ has 
the properties: 

\rm (i) \it there exists a unique $\xi_0\in\Si(1)$ with 
$\langle m,\xi_0\rangle>0$,

\rm (ii) \it there is a cone $\si_0\in\Si$ included in 
$\cal H_m$. 

\nit Then the cone ${\mbb R}_{\geq 0}\xi_0+\si_0$ belongs to $\Si$.}

In {\it loc. cit.} the lemma is stated for $\Si$ complete and 
simplicial and in this case the proof is rather clear. What 
we are going to show here is that the lemma still holds in the 
more general setting when $\Si$ is complete.

\begin{proof}
Let us consider the star of $\si_0$, which is by definition the set 
$$
{\rm Star}(\si_0):=
\left\{
\si\in\Si\mid\si_0\subset\si
\right\}.
$$
Since $\Si$ is complete, 
$$
\Si_0:=\left\{ 
\si/\langle\si_0\rangle\mid\si\in{\rm Star}(\si_0) 
\right\}
$$
is a complete fan in $N_{\mbb R}/\langle\si_0\rangle$. We have used 
the symbol $\langle\si_0\rangle$ to denote the vector space generated 
by the cone $\si_0$. Since $\si_0\in\Ker(m)$, $m$ descends to a linear 
functional $\hat m$ on the quotient
$$
\xymatrix{
N_{\mbb R}\ar[r]^-m\ar[d]& {\mbb R}\\
N_{\mbb R}/\langle\si_0\rangle\ar[ur]_-{\hat m}&
}
$$
and the correponding hyperplane $\cal H_{\hat m}$ breaks 
$N_{\mbb R}/\langle\si_0\rangle$ into two half-spaces. As 
$\Si_0$ is complete, it must exist an element of $\Si_0(1)$ 
situated on the $\{\langle\hat m,\ \cdot\ \rangle>0\}$-side 
of $\cal H_{\hat m}$. But we know from the hypothesis that 
the only element of $\Si(1)$ which is on the 
$\{\langle m,\ \cdot\ \rangle>0\}$-side of $\cal H_m$ is 
$\xi_0$. We deduce that the cone generated by 
$\hat\xi_0:={\rm proj}_{N_{\mbb R}/\langle\si_0\rangle}(\xi_0)$ 
is necessarily an element of $\Si_0(1)$. This means that there 
is a cone $\si'\in {\rm Star}(\si_0)\subset\Si$ such that 
$\si'/\langle\si_0\rangle={\mbb R}_{\geq 0}\hat\xi_0$, and 
consequently there is a $1$-dimensional cone of $\si'$ which 
projects onto ${\mbb R}_{\geq 0}\hat\xi_0$. But, we use the 
unicity of $\xi_0$ again, and we deduce that in fact 
${\mbb R}_{\geq 0}\xi_0$ must be $1$-dimensional cone of $\si'$. 

For the moment we have found a cone $\si'$ of $\Si$ with the 
properties
$$
{\mbb R}_{\geq 0}\xi_0+\si_0\subset\si'\quad
(\xi_0\in\si',\;\si_0\subset\si')
$$ 
and
$$
\si'\subset{\mbb R}_{\geq 0}\xi_0+\langle\si_0\rangle\quad
(\si'/\langle\si_0\rangle={\mbb R}_{\geq 0}\hat\xi_0).
$$
We deduce that $\si'$ is included in the 
$\{\langle m,\ \cdot\ \rangle\geq 0\}$-side of $\cal H_m$, 
and consequently $\si'\cap\cal H_m$ is a face of $\si'$; 
in particular, $\si'\cap\cal H_m$ is an element of $\Si$. 
Using the previous two properties of $\si'$, we find that 
$$
\si_0\subset\si'\cap\cal H_m\subset\langle\si_0\rangle.
$$
As both $\si_0$ and $\si'\cap\cal H_m$ are elements of 
$\Si$, one of the defining properties of a fan implies 
that $\si_0\cap\bigl(\si'\cap\cal H_m\bigr)=\si_0$ is a 
face of $\si'\cap\cal H_m$. This is possible only if 
$\si'\cap\cal H_m=\si_0$ because 
$\si'\cap\cal H_m\subset\langle\si_0\rangle$. We conclude 
now that $\si'(1)={\mbb R}_{\geq 0}\xi\cup \si(1)$ and 
consequently ${\mbb R}_{\geq 0}\xi_0+\si_0=\si'$ is an 
element of $\Si$. 
\end{proof}




\begin{thebibliography}{10}

\bibitem{cox} D.~Cox: {\it The homogeneous coordinate ring of 
a toric variety}, J. Algebr. Geom. {\bf 4} (1995), 17-50
\smallskip

\bibitem{da} V.I.~Danilov: {\it The geometry of toric varieties}, 
Russian Math. Surveys {\bf 33} (1978), 97-154
\smallskip

\bibitem{fu} W.~Fulton: {\it Introduction to toric varieties}, 
Annals of Mathematics Studies no. 131, Princeton University 
Press, 1993
\smallskip

\bibitem{ful} W.~Fulton: {\it Intersection Theory}, Springer-Verlag 
Berlin Heidelberg New York Tokyo 1984
\smallskip

\bibitem{ha} M.~Halic: {\it Higher genus curves on toric varieties}, 
preprint AG/0101100
\smallskip

\bibitem{hus} D.~Husemoller: {\it Fibre Bundles} $2^{\rm nd}$ 
ed. Springer-Verlag 1975
\smallskip

\bibitem{kraft} H.~Kraft: {\it $G$-vector bundles and the 
linearization problem}, Canadian Math. Soc. Conference 
Proceedings {\bf 10} (1989), 111-123
\smallskip

\bibitem{op} T.~Oda, H.~Park: {\it Linear Gale transforms and 
the Gelfand-Kapranov-Zelevinskij decompositions}, Toh\^oku 
Math. J. {\bf 43} (1991), 375-399
\smallskip

\bibitem{ot} Ch.~Okonek, A.~Teleman: {\it Gauge theoretical 
Gromov-Witten invariants for toric varieties. Comparing 
virtual fundamental classes}, preprint

\bibitem{sp} E.~Spanier: {\it Algebraic Topology}, McGraw-Hill 
Inc. 1966
\end{thebibliography}
\end{document}